\documentclass[11pt,a4paper,reqno]{amsart}
\usepackage{amsmath} 
\usepackage{amssymb} 
\usepackage{amscd} 
\usepackage{amsxtra}

\newtheorem{theorem}{Theorem}[section]
\newtheorem{proposition}[theorem]{Proposition}
\newtheorem{corollary}[theorem]{Corollary}
\newtheorem{lemma}[theorem]{Lemma}

\newtheorem{definition}[theorem]{Definition}

\newtheorem{remark}[theorem]{Remark}
\newtheorem{assumption}[theorem]{Assumption}

\makeatletter
\def\le{\mathrel{\mathpalette\gl@align<}}
\def\ge{\mathrel{\mathpalette\gl@align>}}
\def\gl@align#1#2{\lower.6ex\vbox
{\baselineskip\z@skip\lineskip\z@
    \ialign{$\m@th#1\hfil##\hfil$\crcr#2\crcr=\crcr}}}
\makeatother

\newcommand{\C}{{\mathbb C}}

\newcommand{\R}{{\mathbb R}}
\newcommand{\Z}{{\mathbb Z}}

\newcommand{\F}{{\mathbb F}}
\newcommand{\Proj}{{\mathbb P}}

\newcommand{\id}{\operatorname{id}}
\newcommand{\End}{\operatorname{End}}
\newcommand{\Aut}{\operatorname{Aut}}
\newcommand{\ad}{\operatorname{ad}}

\newcommand{\mV}{{\mathcal{V}}}
\newcommand{\mW}{{\mathcal{W}}}

\newcommand{\Lie}{\operatorname{Lie}}
\newcommand{\hq}{\hat{q}}
\newcommand{\Qp}{\boldsymbol{p}}
\newcommand{\Qq}{q}
\newcommand{\hQp}{\hat{\boldsymbol{p}}}
\newcommand{\hQq}{\hat{q}}
\newcommand{\hx}{\hat{x}}
\newcommand{\hy}{\hat{y}}
\newcommand{\Eul}{\mathcal{E}}
\newcommand{\adE}{{\operatorname{ad}}(\mathcal{E})}
\newcommand{\hPhi}{\hat{\Phi}}
\newcommand{\hOmega}{\hat{\Omega}}
\newcommand{\hXi}{\hat{\Xi}}
\newcommand{\ovE}{\overline{E}}
\newcommand{\ovD}{\overline{D}}
\newcommand{\ovPhi}{\overline{\Phi}}
\newcommand{\ovXi}{\overline{\Xi}}
\newcommand{\ovOmega}{\overline{\Omega}}
\newcommand{\barop}{\overline{\phantom{A}}}
\newcommand{\Pic}{\operatorname{Pic}}
\newcommand{\Hom}{\operatorname{Hom}}
\newcommand{\Hol}{\operatorname{Hol}}
\newcommand{\cok}{\operatorname{Coker}}
\newcommand{\rank}{\operatorname{rank}}

\newcommand{\XSig}{{X_{\Sigma}}}

\newcommand{\virtualdim}{\operatorname{virt.dim}}

\def\ov#1{\overline{#1}}

\def\defequal{:=}
\def\parti#1{{q^#1\frac{\partial}{\partial q^#1}}}
\def\hparti#1{{\hat{q}^#1\frac{\partial}{\partial\hat{q}^#1}}}
\def\Qparti#1{{\boldsymbol{q}^#1\frac{\partial}{\partial\boldsymbol{q}^#1}}}
\def\lparfrac#1#2{{\frac{\partial #1}{\partial\log #2}}}
\def\parfrac#1#2{{\frac{\partial #1}{\partial #2}}}
\def\pair#1#2{\langle #1,#2\rangle}
\def\pairV#1#2{\langle #1,#2\rangle^{\mathcal{V}}}

\begin{document}

\title{Quantum $D$-modules 
and Equivariant Floer Theory for Free Loop Spaces}
\author{Hiroshi Iritani} 

\begin{abstract}
The objective of this paper is to clarify 
the relationships between the quantum $D$-module 
and equivariant Floer theory. 
Equivariant Floer theory was introduced 
by Givental in his paper ``Homological Geometry''. 
He conjectured that the quantum $D$-module of a symplectic manifold is 
isomorphic to the equivariant Floer cohomology for 
the universal cover of the free loop space. 
First, motivated by the work of Guest, 
we formulate the notion of ``abstract quantum $D$-module''
which generalizes the $D$-module defined 
by the small quantum cohomology algebra. 
Second, we define the equivariant Floer cohomology of 
toric complete intersections rigorously as a $D$-module, 
using Givental's model. 
This is shown to satisfy the axioms of abstract quantum $D$-module. 
By Givental's mirror theorem \cite{givental-mirrorthm-toric}, 
it follows that equivariant Floer cohomology defined here is 
isomorphic to the quantum cohomology $D$-module.  

\end{abstract}

\maketitle 

\section{Introduction} 
It is known that the quantum cohomology ring 
is isomorphic to Floer cohomology with 
the `pair of pants' product
\cite{piunikhin-salamon-schwartz,ruan-tian}.
Givental proposed the equivariant version of this isomorphism, 
and conjectured that the quantum $D$-module is 
isomorphic to $S^1$-equivariant Floer cohomology
\cite{givental-homologicalgeometry}. 
The quantum $D$-module is a quantum cohomology ring endowed 
with a $D$-module structure. 
On the other hand, 
the $S^1$-equivariant Floer theory 
is the semi-infinite cohomology theory of 
the universal cover $\widetilde{LM}$ of free loop spaces, 
where $S^1$ acts on $\widetilde{LM}$ by rotating loops. 
It (conjecturally) has a $D$-module structure. 
If $M$ is simply-connected, 
any 2-dimensional cohomology class on $M$ can be lifted 
to an equivariant class on $\widetilde{LM}$. 
Equivariant Floer cohomology is considered to 
become a $D$-module 
by multiplication of 2-dimensional equivariant classes and 
pull-back by covering transformations.   
In the paper ``Homological Geometry'', 
Givental constructed an explicit model for $\widetilde{L\Proj_n}$, 
and by a formal computation, 
obtained a generating function of the quantum $D$-module, 
which is called the ``$J$-function''. 
Equivariant Floer theory is, however, not defined rigorously, 
therefore his method has not been justified yet. 
Nevertheless, 
his formal computation of the $J$-function 
gives the correct answer.   
This was proved by another method using the localization
over the moduli space of stable maps by Givental  
\cite{givental-mirrorthm-projective,givental-mirrorthm-toric}
and also by Lian, Liu and Yau 
\cite{lian-liu-yau1,lian-liu-yau2,lian-liu-yau3}. 

Our motivation is to justify Givental's method in the original form.
For that purpose, 
we must first define equivariant Floer theory rigorously. 
In this paper, we construct it 
for complete intersections in toric manifolds 
explicitly by using Givental's model $L$ for $\widetilde{LM}$. 
Floer cohomology is constructed as 
the inductive limit of cohomology algebras.  
This definition makes the meaning of the word 
``semi-infinite dimension'' clear. 
For example, each element in Floer cohomology 
can be expressed as a differential form of infinite degree. 
For non-equivariant Floer cohomology, 
the same construction was done by Cohen-Jones-Segal \cite{cohen-jones-segal} 
in the case of $\Proj^n$. 
The present paper generalizes their construction to intersections in
toric varieties. 
Equivariant Floer cohomology also gives a mathematical model 
for the genus 0 part of Witten's gauged linear sigma models \cite{witten}.

We show that our equivariant Floer cohomology  
has the structure of an {\it{abstract}} quantum $D$-module. 
An abstract quantum $D$-module is, roughly, a  
sheaf over the infinitesimal neighborhood of $0$ in 
$\C^r$ endowed with a flat connection 
with a parameter $\hbar$. 
The original quantum $D$-module is a trivial sheaf over   
$H^2(X,\C^*)$ with fiber $H^*(X,\C)$ 
endowed with a dual Givental connection.  
We choose a partial compactification 
$H^2(X,\C^*)\hookrightarrow \C^r$ by adding a 
``large radius'' limit point $0\in \C^r$. 
The main difference from the original one is that 
we do not a priori fix a frame of the sheaf and 
coordinates $q=(q^1,\dots,q^r)$ of the base space $\C^r$. 
In spite of the indeterminacy, 
we show that these choices are canonically determined 
by the condition that the connection matrices 
represented in the frame are $\hbar$-independent. 
Original quantum $D$-modules have $\hbar$-independent 
connection matrices, 
therefore, from the beginning, they are endowed with 
canonical frames and coordinates.  
In particular, it turns out that 
the abstract $D$-module structure  
determines a quantum deformation of the cup product. 
Guest proved the existence of a canonical frame 
in the neighborhood of a point $q\neq 0$ 
in \cite{guest-Dmodule} using 
the Birkhoff factorization of a loop group. 
We prove the existence and uniqueness of the canonical frame 
in the neighborhood of $q=0$. 
The idea to transform a system of differential equations 
to a ``normal'' form goes back to Birkhoff. 
In singularity theory, the Birkhoff factorization 
appears in connection with 
the construction of primitive forms \cite{saito}. 
In the context of mirror symmetry, 
it appears in the work of Barannikov \cite{barannikov} 
and Coates and Givental \cite{coates-givental}. 
The present paper makes the role of the Birkhoff factorization 
explicit. 


Finally, we show that 
our equivariant Floer cohomology
is isomorphic to the quantum $D$-module for 
nef toric complete intersections
using Givental's mirror theorem in 
\cite{givental-mirrorthm-toric,lian-liu-yau3}. 
This follows from the coincidence of 
the two $J$-functions, one from equivariant Floer cohomology and 
the other from quantum cohomology. 
We will treat the case where the first Chern class is not nef 
in \cite{iritanigen}.  

Strictly speaking, we regard a toric complete intersection 
as a superspace $(\XSig,\mV)$ --- a pair consisting of a
toric manifold  $\XSig$ and a sum $\mV$ of nef line bundles. 
We will introduce the notion of quantum products and  
quantum $D$-module for superspaces 
in Section \ref{sect:quantumcohomologyandDmodule}. 
It takes into consideration only quantum multiplications of classes from 
the ambient toric manifolds. 
Therefore, the quantum $D$-module under consideration is 
for the superspace $(\XSig,\mV)$ rather than for the 
complete intersection itself. 

Another important point in this paper is   
a pairing between equivariant Floer homology and cohomology. 
It is defined as the infinite dimensional 
Poincar\'{e} pairing between semi-infinite cycle and cocycles.  
Moreover, this geometrical pairing coincides with the 
pairing defined for abstract quantum $D$-modules.  

The paper is organized as follows. 
In Section 2, we explain quantum cohomology and $D$-modules 
for convex superspaces.  
In Section 3, we formulate the abstract quantum $D$-module. 
We prove the existence and 
uniqueness of the canonical frame.  
In Section 4, we construct the equivariant Floer cohomology for 
toric complete intersections. 
We prove that the equivariant Floer cohomology is 
isomorphic to the quantum $D$-module for nef toric complete intersections.  
In Section 5, we illustrate the general theory by 
examples using the Hirzeburch surface. 

We remark that in this paper, 
we only deal with the small quantum cohomology, not the large one. 
We also restrict our attention to manifolds whose 
cohomology rings are generated by two dimensional classes. 
A generalization to the large quantum cohomology 
will also be treated in \cite{iritanigen}.

\noindent
{\bf Acknowledgments} Thanks are due to 
Professor Hiraku Nakajima, my advisor, 
for many helpful comments and discussions 
and also to Professor Toshiaki Maeno, Professor Kenji Fukaya 
for valuable discussions on Homological Geometry.  
The author expresses gratitude to Naohiro Sakamoto.  
The idea of taking the inductive limit 
was suggested by him. 
The author is also grateful to Professor Kyoji Saito and Atsushi Takahashi 
for their helpful comments.  
The author expresses deep gratitude to Professor Martin Guest 
for reading the manuscript very carefully and 
finding many English mistakes and for valuable comments.  
This research was partially supported by 
Grant-in-Aid for JSPS Fellows  
and Scientific Research (15-5482).

\section{Quantum Cohomology and $D$-modules}
\label{sect:quantumcohomologyandDmodule}

In this section we introduce the notion of  
quantum $D$-module for convex supermanifold. 
Here, we mean by the word supermanifold 
a pair $(M,\mV)$ of a manifold $M$ 
and a vector bundle $\mV$ on it, 
following Givental's terminology in 
\cite{givental-mirrorthm-toric}. 
The fundamental class of the supermanifold is 
by definition the cap product of the fundamental class 
of $M$ with the Euler class $Euler(\mV)$ of $\mV$.  

Let $M$ be a smooth projective manifold whose total cohomology 
ring is generated by two dimensional cohomology classes and 
$\mV$ be a holomorphic vector bundle over $M$. 
We assume that $\mV$ is convex, i.e. for any rational curve $C$ and 
holomorphic map $f\colon C\rightarrow M$, we have 
$H^1(C,f^*(\mV))=0$. 
The first Chern class of the supermanifold is 
defined by $c_1(M/\mV)\defequal c_1(TM)-c_1(\mV)$. 
Later, we assume that $c_1(M/\mV)$ is a nef class, 
but this assumption is unnecessary in this section. 
Let $\{T_j\}_{j=0}^{s}$ be a basis of total cohomology ring $H^*(M,\C)$, 
and $\{T^j\}_{j=0}^{s}$ be its dual basis with 
respect to the Poincar\'{e} pairing i.e. 
$\int_{M} T_i\cup T^j=\delta_{i}^j$. 
We assume that $T_0=1$ and $T_s$ is the Poincar\'{e} dual 
of a point. 

Let $\Lambda\subset H_2(M,\Z)_{\text{free}}$ 
be the semigroup generated by classes of effective curves, 
where $H_2(M,\Z)_{\text{free}}$ 
means $H_2(M,\Z)$ modulo its torsion. 
For the superspace $(M,\mV)$, we have a quantum product structure 
on the module $H^*(M,\C)\otimes\C[\![\Lambda]\!]$, 
which was first introduced in \cite{givental-mirrorthm-projective} 
and explicitly written down in \cite{pandharipande-aftergivental}. 
Let $\ov{M}_{0,n}(M,d)$ 
be the moduli space of genus zero,  
degree $d$ stable maps to $M$ with $n$ marked points. 
We have a forgetful map 
$\pi_{i}\colon \ov{M}_{0,n+1}(M,d)\rightarrow \ov{M}_{0,n}(M,d)$ 
which forgets the $i$-th marked point 
and an evaluation map 
$e_i\colon \ov{M}_{0,n+1}(M,d)\rightarrow M$ 
which evaluates at the $i$-th marked point. 
We define an orbibundle $\mV_d$ over the moduli space 
$\ov{M}_{0,n}(M,d)$ as 
\[\mV_d\defequal {\pi_{n+1}}_{*}e_{n+1}^*(\mV). \]  
The fiber of $\mV_d$ at the stable map $f\colon C\rightarrow M$ 
is isomorphic to $H^0(C,f^*(\mV))$. 
Therefore by the convexity of $\mV$, 
the rank of $\mV_d$ is equal to $\pair{c_1(\mV)}{d}+1$. 
Define a subbundle $\mV'_{d,i}$ by the exact sequence 
\[
\begin{CD}
0  @>>>  \mV'_{d,i} @>>> \mV_d  @>>> e_i^*(\mV) @>>> 0, 
\end{CD}
\] 
where the last map is an evaluation map 
of sections in $H^0(C,f^*(\mV))$ at the $i$-th marked point. 
We introduce a fiber-wise  $S^1$ action on the bundles 
$\mV$, $\mV_d$ and $\mV'_{d,i}$ i.e. $S^1$ acts on each fiber 
by scalar multiplication and trivially on the base. 
Let $\lambda$ be a generator of the equivariant cohomology 
of a point with respect to this $S^1$ action.  
Define the small quantum product $*$ by the formula 
\begin{equation}
\label{eq:quantumproduct} 
T_i*T_j=\sum_{\substack{d \in\Lambda \\ k=0,\dots,s}} q^d T^k
  \int_{[\ov{M}_{0,3}(M,d)]^{\text{virt}}}
  e_1^*(T_i)\cup e_2^*(T_j)\cup e_3^*(T_k)\cup Euler_{S^1}(\mV'_{d,3}),
\end{equation}
where $q^d$ denotes an element in the group algebra $\C[\Lambda]$. 
This product structure $*$ is extended linearly over 
$\C[\lambda][\![\Lambda]\!]$ 
and defines a commutative, associative algebra structure 
on the module 
$H^*(M,\C)\otimes \C[\lambda][\![\Lambda]\!]$.  
We call this algebra the {\bf (equivariant) small quantum cohomology} 
of the superspace $(M,\mV)$ and denote it by 
$QH^*_{S^1}(M/\mV)$.  
Similarly we can define its non-equivariant counterpart $QH^*(M/\mV)$. 
The product $*$ is a $q$-deformation of the usual cup product. 
If we assign degree $2\pair{c_1(M/\mV)}{d}$ to an element $q^d$ 
in $\C[\Lambda]$, then 
$QH^*_{S^1}(M/\mV)$ and $QH^*(M/\mV)$ become graded rings.  
Define a Poincar\'{e} pairing for the superspace by 
\[\pairV{T_i}{T_j}=\int_{M}T_i\cup T_j\cup Euler_{S^1}(\mV).\]
Then we have the following Frobenius property. 
\begin{equation}
\label{eq:Frobenius}
  \pairV{T_i}{T_j*T_k}=\pairV{T_i*T_j}{T_k} \quad 
  \text{for all } i,j,k. 
\end{equation}
When we have $\mV=0$, the above definition of quantum cohomology 
of $(M,\mV)$ coincides with that of $M$ itself.  
The quantum product $*$ of superspaces is less familiar 
than that of projective manifolds.  
In fact, it is closely related to the quantum cohomology 
of the zero-locus $N$ of a transverse section of $\mV$. 
Let $i\colon N\rightarrow M$ be the inclusion. 
By \cite[Proposition 4]{pandharipande-aftergivental}, we have 
a ring homomorphism 
\[i^*\colon QH^*(M/\mV)\longrightarrow QH^*(N), \]
provided that the tangent bundle $TM$ is convex 
and $H_2(M,\Z)\cong H_2(N,\Z)$.   
Therefore we can say that quantum cohomology for superspaces 
detects quantum products of classes in $N$ which come from 
the ambient space $M$. 

For convenience, 
we choose coordinates of $\C[\Lambda]$ and 
modify the above definition of small quantum cohomology slightly.
The set of nef classes in $H^2(M,\Z)_{\text{free}}$ forms a semigroup
and its convex hull in $H^2(M,\R)$ forms a cone $\sigma$. 
The interior of the cone $\sigma$ is non-empty because it 
contains a class represented by the K\"{a}hler form. 
We choose nef classes 
$p_1,\dots,p_r$ in $H^2(M,\Z)_{\text{free}}$ 
which form a $\Z$-basis of $H^2(M,\Z)_{\text{free}}$. 
When $c_1(M/\mV)$ is nef, we choose $p_1,\dots,p_r$ so that  
$c_1(M/\mV)$ is contained in the cone generated by $p_1,\dots,p_r$. 
Let $q^1,\dots,q^r$ in $H_2(M,\Z)_{\text{free}}$ be a dual basis of 
$p_1,\dots,p_r$ 
satisfying $\langle p_a,q^b\rangle =\delta_a^b$.  
For an element $d\in \Lambda$, $q^d\in \C[\Lambda]$ is identified with 
the element 
\[(q^1)^{\langle p_1,d\rangle}(q^2)^{\langle p_2, d\rangle}
\dots (q^r)^{\langle p_r,d\rangle}\quad 
\text{ in }\C[q^1,\dots,q^r].\]
Therefore, we can embed $\C[\Lambda]$ into 
the polynomial ring $\C[q^1,\dots,q^r]$.
Note that $\langle p_a,d\rangle$ is non-negative because $p_a$  
is nef. 
Through the above inclusion, 
the degree of $q^a$ is identified with $2\langle c_1(M),q^a\rangle$. 
{\it Note that $\deg q^a$ is non-negative if $c_1(M/\mV)$ is nef}. 
In general, the degree of $q^a$ can be negative. 
We frequently regard $q^a$ as a complex coordinate of 
$H^2(M,\C^*)$ as follows.  
\begin{eqnarray*}
q^a\colon H^2(M,\C^*)=H^2(M,\C/2\pi \sqrt{-1}\Z)
\ni\alpha\longmapsto\exp(\langle \alpha,q^a\rangle)
\in \C
\end{eqnarray*}
After choosing a set of coordinates $q^1,\dots,q^r$, we replace 
$\C[\Lambda]$ with $\C[q^1,\dots,q^r]$, and set 
$QH^*_{S^1}(M/\mV)
\defequal H^*(M,\C)\otimes\C[\lambda][\![q^1,\dots,q^r]\!]$, 
where $\C[\lambda][\![q^1,\dots,q^r]\!]$ denotes 
the ring of formal power series in $q^1,\dots,q^r$
with coefficients in $\C[\lambda]$. 
For simplicity, we denote it by $\C[\lambda][\![q]\!]$.   
The product $*$ is defined by the same formula 
$(\ref{eq:quantumproduct})$.  

Next we define the quantum $D$-module of the superspace $(M,\mV)$. 
The {\bf dual Givental connection} 
is a connection defined on the trivial bundle 
$H^*(M,\C)\times H^2(M,\C^*)\rightarrow H^2(M,\C^*)$
by 
\[\nabla^{\hbar}\defequal \hbar d+\sum_{a=1}^r (p_a*)\frac{dq^a}{q^a}, 
\]
where $\hbar$ is a degree two parameter. 
Here, we follow the convention that 
the Givental connection means $-\nabla^{-\hbar}$
and the `dual' Givental connection means $\nabla^{\hbar}$ 
\cite[p.311,p.321]{cox-katz}.  
This is a formal connection. 
Algebraically, it defines a map 
\[\nabla^{\hbar}\colon H^*(M)\otimes\C[\hbar,\lambda][\![q]\!]
 \longrightarrow \sum_a H^*(M)\otimes \C[\hbar,\lambda][\![q]\!] 
 \frac{dq^a}{q^a}. 
\]
Put $\nabla^{\hbar}_a\defequal \nabla^{\hbar}_{\parti{a}}
=\hbar\parti{a}+p_a*$. 
The connection $\nabla^{\hbar}$ is flat.
This fact may be verified from the existence of a fundamental solution 
matrix $L$ which satisfies $\nabla^{\hbar}L=0$. 
A fundamental solution $L$ is 
given in \cite[Equation (25)]{pandharipande-aftergivental} 
explicitly as follows. 
\begin{equation}
\label{eq:non-abstractL} 
L(T_i)\defequal e^{-p\log q/\hbar}T_i-
       \sum_{\substack{d\in \Lambda\setminus\{0\}\\ k=0,\dots,s}}
       q^d T^k\int_{[\ov{M}_{0,2}(M,d)]^{\text{virt}}}
       \frac{e_1^*(e^{-p\log q/\hbar}T_i)}{\hbar+\psi_1}e_2^*(T_k)  
       Euler_{S^1}(\mV'_{d,2}),
\end{equation}
where $\psi_1$ denotes the first Chern class of the first cotangent line 
bundle on $\ov{M}_{0,2}(M,d)$ and 
$p\log q\defequal \sum_{a=1}^{r}p_a\log q^a$. 
The expression $(\hbar+\psi_1)^{-1}$ in the integrand may be 
written more explicitly as  
\[\sum_{j=1}^{\infty} \hbar^{-j-1}(-\psi_1)^j.\]
The sign of $\hbar$ in (\ref{eq:non-abstractL})
is opposite to that in \cite[Equation (25)]{pandharipande-aftergivental}
because we use the dual Givental connection 
instead of the Givental connection.   
Let $D$ denote a Heisenberg algebra 
\[D\defequal 
  \C[\hbar][\![\Qq^1,\dots,\Qq^r]\!][\Qp_1,\dots,\Qp_r]  
\quad \deg \Qq^i=2\langle c_1(M/\mV),q^i\rangle,\;
\deg \Qp_i=\deg\hbar=2,\]
whose generators satisfy the following commutation relations.
\begin{gather*}
[\Qp_a,\Qq^b]=\hbar\delta_a^b\Qq^b,\;
[\Qp_a,\Qp_b]=[\Qq^a,\Qq^b]=0, \\
[\Qp_a,f(\Qq^1,\dots,\Qq^r)]=\hbar\Qparti{a}f(\Qq^1,\dots,\Qq^r)
\text{ for }f\in \C[\![\Qq^1\dots,\Qq^r]\!].
\end{gather*}
We define an action of $D$ on the module 
$E_{S^1}\defequal H^*(M,\C)\otimes \C[\hbar,\lambda][\![q]\!]$ by 
\[\Qq^a\mapsto \text{multiplication by }q^a,\qquad 
\Qp_a\mapsto \nabla^\hbar_a.\]
We call this $D$-module $E_{S^1}$ 
the {\bf (equivariant) small quantum $D$-module} 
arising from the quantum cohomology $QH^*_{S^1}(M/\mV)$. 
Similarly, we can define its non-equivariant counterpart 
$E=H^*(M)\otimes\C[\hbar][\![q]\!]$.  

Next we introduce the $J$-function of the quantum $D$-module. 
The fundamental solution matrix $L$ can be decomposed as 
$L=S\circ e^{-p\log q/\hbar}$, where 
\begin{equation}
\label{eq:non-abstractS}
S(T_i)\defequal T_i-
       \sum_{\substack{d\in \Lambda\setminus\{0\}\\ k=0,\dots,s}}
       q^d T^k\int_{[\ov{M}_{0,2}(M,d)]^{\text{virt}}}
       \frac{e_1^*(T_i)}{\hbar+\psi_1}e_2^*(T_k)  
       Euler_{S^1}(\mV'_{d,2}).
\end{equation}
Note that $S$ is an element of 
$\End(H^*(M))\otimes \C[\hbar^{-1},\lambda][\![q]\!]$. 
Define the $J$-function by 
\[J(q,\hbar)\defequal L^{-1}(1)=e^{p\log q/\hbar}S^{-1}(1).\]
The function $J(q,\hbar)$ is a cohomology valued formal function.
Next, we show that this $J$-function  
coincides with that given in 
\cite{givental-mirrorthm-projective,
givental-mirrorthm-toric, cox-katz}.
 
\begin{lemma}
\label{lem:non-abstractdiffeqforS} 
Let $\mu$ be a constant matrix defined by $\mu(T_i)=(\deg T_i)T_i$. 
The matrix function $S$ satisfies the following differential equations. 
\begin{gather*}
\hbar\parti{a}S+(p_a*)\circ S-S\circ (p_a\cup)=0. \\ 
(2\hbar\frac{\partial}{\partial\hbar}
        +\sum_{a=1}^r (\deg q^a)\parti{a})S+[\mu,S]=0. 
\end{gather*}
\end{lemma} 

The first equation follows from $\nabla^{\hbar}L=0$. 
The second one follows from the fact that the matrix $S$ 
preserves the degree.   
Using the above equations, we can establish the 
unitarity of $S$ which is stated in \cite[p.117]{givental-elliptic}.  

\begin{lemma}
\label{lem:Sisunitary}
For $\alpha,\beta\in H^*(M)$, we have 
\[\pairV{S(\alpha)(q,-\hbar)}{S(\beta)(q,\hbar)}=\pairV{\alpha}{\beta}.\] 
\end{lemma} 
\begin{proof} 
Set the left hand side equal to $F(\alpha,\beta)$. 
First, $F(\alpha,\beta)|_{q=0}=\pairV{\alpha}{\beta}$ is clear. 
By using the differential equation satisfied by $S$ and the 
Frobenius property, we have 
\[\hbar\parti{a}F(\alpha,\beta)=
  -F(p_a\cup \alpha,\beta)+F(\alpha,p_a\cup\beta).
\] 
Because the matrix $(p_a\cup)$ is nilpotent, there exists some $n$ 
such that $(\hbar\parti{a})^nF(\alpha,\beta)=0$ holds.  
Therefore the lemma holds.  
\end{proof}

The $J$-function can be explicitly calculated as follows.  
\begin{proposition}
\label{prop:explicitJ}
The inverse matrix $S^{-1}$ has the following explicit formula. 
\[S^{-1}(T_i)=T_i+\sum_{\substack{d\in\Lambda\setminus \{0\} \\
                            j=0,\dots,s}}q^d T^j
       \int_{[\ov{M}_{0,2}(M,d)]^{\text{\rm virt}}}
         \frac{e_1^*(T_j)}{\hbar-\psi_1}e_2^*(T_i)Euler_{S^1}(\mV'_{d,1}).
\]
In particular, we have 
\[J(q,\hbar)=e^{p\log q/\hbar}\Big
           (1+\sum_{\substack{d\in\Lambda\setminus \{0\} \\ 
                            j=0,\dots,s}}q^d T^j
      \int_{[\ov{M}_{0,2}(M,d)]^{\text{\rm virt}}}
        \frac{e_1^*(T_j)}{\hbar-\psi_1}Euler_{S^1}(\mV'_{d,1})\Big).
\]
\end{proposition}
\begin{proof} By the unitarity of $S$ and the non-degeneracy of 
the equivariant pairing $\pairV{\cdot}{\cdot}$, 
it suffices to calculate $\pairV{S(\alpha)(q,-\hbar)}{\beta}$. 
Using the equalities 
$\beta\cup Euler_{S^1}(\mV)=\sum_{j=0}^s\pairV{T^j}{\beta}T_j$
and 
$e_i^*(Euler_{S^1}(\mV))\cup Euler_{S^1}(\mV'_{d,i})=Euler_{S^1}(\mV_d)$, 
we have 
\begin{align*}
\pairV{S(\alpha)(q,-\hbar)}{\beta}&=\pairV{\alpha}{\beta}
  +\sum_{\substack{d\in\Lambda\setminus \{0\}\\ j=0,\dots,s}} 
    q^d \pairV{T^j}{\beta} \int_{[\ov{M}_{0,2}(M,d)]^{\text{virt}}}
    \frac{e_1^*(\alpha)}{\hbar-\psi_1}e_2^*(T_j) Euler_{S^1}(\mV'_{d,2}) \\
&=\pairV{\alpha}{\beta}
   +\sum_{d\in\Lambda\setminus \{0\}} 
    q^d \int_{[\ov{M}_{0,2}(M,d)]^{\text{virt}}}
    \frac{e_1^*(\alpha)}{\hbar-\psi_1}e_2^*(\beta\cup 
    Euler_{S^1}(\mV)) Euler_{S^1}(\mV'_{d,2})  \\
&=\pairV{\alpha}{\beta}
   +\sum_{d\in\Lambda\setminus \{0\}} 
    q^d \int_{[\ov{M}_{0,2}(M,d)]^{\text{virt}}}
    \frac{e_1^*(\alpha\cup Euler_{S^1}(\mV))}
         {\hbar-\psi_1}e_2^*(\beta) Euler_{S^1}(\mV'_{d,1}) \\
&=\pairV{\alpha}{\beta}
   +\sum_{\substack{d\in\Lambda\setminus \{0\}\\ j=0,\dots,s}} 
    q^d \pairV{\alpha}{T^j}\int_{[\ov{M}_{0,2}(M,d)]^{\text{virt}}}
    \frac{e_1^*(T_j)}{\hbar-\psi_1}e_2^*(\beta) Euler_{S^1}(\mV'_{d,1}). 
\end{align*}
This gives the formula for $S^{-1}$. 
\end{proof}

The formula for the $J$-function coincides with \cite[p.381]{cox-katz} 
up to multiplication by the Euler class.
Note that functions $S$, $L$ and $J$ can be defined also in the 
non-equivariant theory. 
The reason for introducing the $S^1$-equivariance 
is to make the pairing $\pairV{\cdot}{\cdot}$ non-degenerate.  

The $J$-function plays an important role in Givental's theory. 
We shall show next that the $J$-function is a generator of the 
quantum $D$-module. 
Note that the following theorem holds true without 
the nef assumption for $c_1(M/\mV)$. 

\begin{theorem} 
\label{thm:Jisgen}
Let $M$ be a smooth projective manifold 
and $\mV$ be a convex vector bundle on $M$. 
Assume that the total cohomology ring of $M$ 
is generated by the second cohomology group.   
Then the $J$-function of the supermanifold $(M,\mV)$ is a 
generator of its small quantum $D$-module 
$E=H^*(M,\C)\otimes\C[\hbar][\![q]\!]$. 
More precisely, we have a $D$-module isomorphism 
$D/I\cong E$, 
where $I$ is the left ideal consisting of elements 
$f(\Qq,\Qp,\hbar)$ in $D$ satisfying
\begin{equation}
\label{eq:killJ}
f(q,\hbar\parti{{}},\hbar)\cdot J=0.
\end{equation}
\end{theorem} 

\begin{corollary}
Let $(M,\mV)$ be a supermanifold which satisfies 
the assumptions in Theorem \ref{thm:Jisgen}.   
The isomorphism type of the quantum cohomology ring 
of the superspace $(M,\mV)$ is completely determined by its 
$J$-function. More precisely, we have an isomorphism
\[\C[\![q^1,\dots,q^r]\!][p_1,\dots,p_r]
/\lim_{\hbar\rightarrow 0}I\cong QH^*(M/\mV),\] 
where 
$\lim_{\hbar\rightarrow 0}I \defequal 
\{f(q,p,0)\;|\; f(\Qq,\Qp,\hbar)\in I\}$ 
and $q^a$ and $p_a$ are commuting variables corresponding to 
$\Qq^a$ and $\Qp_a$.  
(Note that the non-commutativity vanishes when $\hbar\to 0$.) 
\end{corollary} 

Later, we will see that the $J$-function also determines 
small quantum products, i.e. a matrix representation 
of the quantum product $p_a*$. 
Moreover, by the reconstruction theorem of Kontsevich and Manin, 
we see that the $J$-function determines all genus 0 $n$-points   
Gromov-Witten invariants when $\mV=0$.

\noindent
{\it Proof of Theorem \ref{thm:Jisgen}}. 
First note that $L$ is an intertwining operator of two $D$-module
structures, i.e. $\nabla^{\hbar}\circ L=L\circ \hbar d$. 
From the definition $J=L^{-1}(1)$, 
we have a $D$-module injective homomorphism $D/I\rightarrow E$ 
which sends $f(\Qq,\Qp,\hbar)$ in $D$ 
to $f(q,\nabla^{\hbar},\hbar)\cdot 1$ in $E$.  
It suffices to show that it is surjective. 
We show that for $\alpha\in H^*(M,\C)$, 
there exists $f(\Qq,\Qp,\hbar)\in D$ 
such that $\alpha=f(q,\nabla^{\hbar},\hbar)\cdot 1$. 
Take a polynomial $P(x_1,\dots,x_r)$ such that 
$\alpha=P(p_1,\dots,p_r)$ holds in the cohomology ring. 
Noting the relation 
$\nabla^{\hbar}_a(g(q)\beta)\equiv g(0)p_a\cup\beta 
\mod \langle q^1,\dots,q^r\rangle $ 
for $\beta\in H^*(M)$ and $g(q)\in \C[\![q]\!]$, 
we can write  
\[P(\nabla^{\hbar}_1,\dots,\nabla^{\hbar}_r)\cdot 1
    = \alpha -\sum_{i=1}^r q^i \beta_i(\hbar) + 
      \text{\rm higher order terms in $q$}, 
\]
for some $\beta_i(\hbar)\in H^*(M)[\hbar]$. 
Next we take a polynomial $P_i$  
such that $P_i(p_1,\dots,p_r,\hbar)=\beta_i(\hbar)$ holds. 
We have 
\[\{P(\nabla^{\hbar})
  +\sum_{i=1}^r q^i P_i(\nabla^{\hbar},\hbar)
  \}\cdot 1
  =\alpha - \sum_{i,j=1}^r q^i q^j \beta_{ij}(\hbar)
          +\text{\rm higher order terms in $q$}.
\]
The degree of each polynomial $P(p)$, $P_i(p,\hbar)$ with respect to 
$p_1,\dots,p_r$ is bounded from above by $\dim_{\C} M$. 
Therefore, repeating this process, 
we obtain $f(\Qq,\Qp,\hbar)\in \C[\hbar][\![\Qq]\!][\Qp]$ 
such that $f(q,\nabla^\hbar,\hbar)\cdot 1=\alpha$ holds. 
\qed


\section{Abstract Quantum $D$-modules }
\label{sect:abstractquantumDmodule}

In this section, we formulate the abstract quantum $D$-module 
over the base space $B\defequal \C^r$. 

\subsection{Definitions and Notation}

Let $\{q^1,\dots,q^r,\hbar\}$ be a coordinate system 
of $B\times \C=\C^r\times \C$ centered at $0$. 
Let $\mathcal{O}$ and $\mathcal{O}^\hbar$ denote 
the rings of germs of regular functions at $0$ 
on $B$ and $B\times \C$ respectively. 
We can consider the abstract quantum $D$-module 
both in the formal and convergent categories. 
In the formal category, we take 
\[\mathcal{O}=\C[\![q^1,\dots,q^r]\!], \quad 
  \mathcal{O}^{\hbar}=\C[\hbar][\![q^1,\dots,q^r]\!], 
\]
and in the convergent one, we take $\mathcal{O}^{\hbar}$ 
to be the ring of convergent power series in $\hbar$ and $q^1,\dots, q^r$. 
However, for simplicity, we deal only with formal power series 
and do not discuss their convergence. 
Most of the results hold also in the convergent category under some 
additional conditions. 

Let $D$ be the Heisenberg algebra of the previous section. 
\[D=\mathcal{O}^\hbar[\Qp_1,\dots,\Qp_r]
   =\C[\hbar][\![\Qq^1,\dots,\Qq^r]\!][\Qp_1,\dots, \Qp_r].  
\]
We begin with the definition of generalized coordinates of $D$. 
$\mathcal{O}$ can be considered as a subring of $D$. 
Define two left $\mathcal{O}$-submodules 
$\mathfrak{m}$, $\mathfrak{m}'$ of $D$ as 
$\mathfrak{m}=\sum_{a=1}^r\mathcal{O}q^a$  
and $\mathfrak{m}'=\mathfrak{m}+\sum_{a=1}^r\mathcal{O}\Qp_a$. 
Here, $\mathfrak{m}$ is the ideal of 0 in $\mathcal{O}$, 
and $\mathfrak{m}'$ is a non-commutative analogue of the ideal of $0$ 
(but is not an ideal of $D$). 
A set of elements 
$(\hQq^1,\dots,\hQq^r,\hQp_1,\dots,\hQp_r)$ in $D$ is called 
{\bf generalized coordinates} of $D$ if it satisfies the following conditions. 
\begin{enumerate} 
\item  $\hQq^1,\dots,\hQq^r$ generate $\mathfrak{m}$ as 
       an $\mathcal{O}$-module and  
       $\hQq^1,\dots,\hQq^r,\hQp_1,\dots,\hQp_r$ generate $\mathfrak{m}'$ as 
       an $\mathcal{O}$-module. 
\item  $[\hQp_a,\hQq^b]=\hbar \delta_a^{b}\hQq^b$, 
       $[\hQp_a,\hQp_b]=0$. 
\end{enumerate}
For example, $(\Qq^1,\dots,\Qq^r,\Qp_1,\dots,\Qp_r)$ is 
a set of generalized coordinates of $D$. 
We can define the notion of generalized coordinates more concretely,
as the following proposition shows. 

\begin{proposition}
\label{prop:coordinateofD}
Let $(\hQq^1,\dots,\hQq^r,\hQp_1,\dots,\hQp_r)$ 
be a set of generalized coordinates of $D$. 
Then $(\hQq^1,\dots,\hQq^r)$ gives a local coordinate system of $B$ 
centered at $0$ and the $\log$-Jacobi matrix 
\[\frac{\partial\log \hq^a}{\partial\log q^b} 
\defequal\frac{q^b}{\hq^a}\frac{\partial\hq^a}{\partial q^b}\] 
is non-degenerate and regular in the neighborhood of 
$0$. Here, the word ``regular'' means each entry 
is contained in $\mathcal{O}$. 
Moreover, $\hQp_a$ can be written as 
\begin{equation}
\label{eq:coordinatetransformation} 
\hQp_a=\sum_{b=1}^r\lparfrac{\log q^b}{\hq^a}\Qp_b+\hparti{a}F 
\end{equation} 
for some element $F$ in $\mathfrak{m}$. 
\end{proposition} 
We need the following lemma for the proof. 

\begin{lemma}
\label{lem:integrabilitycond}
Given $F_1,\dots,F_r\in \mathfrak{m}$, if we have 
$\hparti{a}F_b=\hparti{b}F_a$, 
there exists an element $F$ in $\mathcal{O}$ such that 
$F_a=\hparti{a}F$. 
\end{lemma} 

\noindent
{\it Proof of Proposition \ref{prop:coordinateofD}}.  
It is clear that $\{\hq^1,\dots,\hq^r\}$ forms 
a set of local coordinates of $B$ centered at $0$. 
$\hQp_a$ can be written in the form, 
\begin{equation*}
\hQp_a=\sum_{b=1}^r G_a^b\Qp_b+F_a,\; 
(G^b_a)\text{ is an invertible matrix with entries in } \mathcal{O},
\; F_a\in \mathfrak{m}
\end{equation*}
By the relations $[\hQp_a,\hQp_b]=0$, 
$[\hQp_a,\hQq^b]=\hbar\delta_a^b\hQq^b$, we have
\begin{equation} 
\label{eq:intheprop:coordinateofD}
G_a^c\parti{c}G_b^d=G_b^c\parti{c}G_a^d,\quad
G_a^c\parti{c}F_b=G_b^c\parti{c}F_a,\quad
G_a^c\parti{c}\hq^b=\delta_a^b\hq^b,
\end{equation}
where we sum over the index $c$. 
The third equation shows that $G_a^b$ is 
the inverse matrix of the $\log$-Jacobi matrix. 
In particular, the $\log$-Jacobi matrix is non-degenerate
and regular and 
\[G_a^b=\lparfrac{\log q^b}{\hq^a}.\] 
Substituting the right hand side for $G_a^c$, $G_b^c$ in the second 
equations of (\ref{eq:intheprop:coordinateofD}), we obtain 
\[\hparti{a}F_b=\hparti{b}F_a.\]
By Lemma \ref{lem:integrabilitycond}, we obtain the proposition. 
\qed 

\begin{remark}
\label{rem:senseoft_0}
Because of the inhomogeneous term $\parti{a}F$ 
in equation (\ref{eq:coordinatetransformation}), 
it seems that we cannot regard $\Qp_a$ simply as a differential operator 
$\hbar\parti{a}$. 
However, introducing a new variable $t^0$ formally and 
considering the following coordinate transformation,
\[t^0=\hat{t}^0+F(\hat{q}),\quad q^a=q^a(\hat{q}),\]
we have the following relations. 
\[\left(\hbar\lparfrac{}{\hq^a}\right)=\sum_{b=1}^r 
 \lparfrac{\log q^b}{\hq^a}\left(\hbar\lparfrac{}{q^b}\right)
 +\lparfrac{F}{\hq^a}\left(\hbar\parfrac{}{t^0}\right).
\]
In this case, 
if we regard $\hbar\parfrac{}{t^0}$ as the identity operator, 
we recover equation (\ref{eq:coordinatetransformation}).  
\end{remark}

We see that any set of coordinates $(\hq^1,\dots,\hq^r)$ of 
$B$ which comes from a set of generalized coordinates 
$(\hQq^1,\dots,\hQq^r,\hQp_1,\dots,\hQp_r)$ of $D$ 
has a special relation to the original coordinates 
$(q^1,\dots,q^r)$ i.e. the $\log$-Jacobi matrix 
is non-degenerate and regular.  
We say that such coordinates of $B$ centered at $0$ 
are {\bf admissible}. 
From now on, we do not fix a generalized coordinate of $D$, 
therefore $(\Qq^1,\dots,\Qq^r,\Qp_1,\dots,\Qp_r)$ 
denotes an arbitrary set of generalized coordinates of $D$ which 
does not necessarily coincide with the 
original one.
We give a necessary and sufficient condition that 
a matrix with entries in $\mathcal{O}$ becomes 
a $\log$-Jacobi matrix of some coordinate transformation. 

\begin{lemma} 
\label{lem:integrabilitycondforlog-Jacobi}
Let $(G^a_b)$ be an invertible $r\times r$ matrix 
with entries in $\mathcal{O}$. 
$(G^a_b)$ is a $\log$-Jacobi matrix of 
some coordinate transformation i.e. 
$G^a_b=\partial\log q^a/\partial \log \hq^b$  
if and only if it satisfies 
\begin{equation}
\label{eq:integrabilitycondforlog-Jacobi}
\hparti{c}G^a_b=\hparti{b}G^a_c,\qquad
G^a_b(\hq=0)=\delta^{\sigma(a)}_{b}
\end{equation}
for a permutation $\sigma\in \mathfrak{S}_r$. 
Given admissible coordinates $(\hq^1,\dots,\hq^r)$ and  
a $\log$-Jacobi matrix $(G^a_b)$, the new set of coordinates 
$(q^1,\dots,q^r)$ is of the form 
\[q^a=\hq^{\sigma(a)}\exp(\delta^a), \quad \delta^a\in\mathcal{O},
\]
and is determined up to constants, 
$q^a\mapsto c^aq^a$.  
\end{lemma} 
\begin{proof}
Set $C^a_b=G^a_b(\hq=0)$.  
Assume that $G^a_b$ satisfies the first equation of 
(\ref{eq:integrabilitycondforlog-Jacobi}).  
Because $G^a_b-C^a_b \in\mathfrak{m}$, 
Lemma \ref{lem:integrabilitycond} shows that  
there exist elements $\delta^a\in\mathcal{O}$
satisfying 
$\hparti{b}\delta^a=G^a_b-C^a_b$. 
Thus, the solution for 
$G^a_b=\partial\log q^a/\partial \log \hq^b$ is given by 
$q^a=\prod_{b=1}^r (\hq^b)^{C^a_b}\exp(\delta^a)$. 
In order for $(q^1,\dots,q^r)$ to be a set of coordinates, 
$(C^a_b)$ must be an element of $SL(r,\Z)$ and 
$C^a_b\ge 0$, ${C^{-1}}^a_b\ge 0$. 
Therefore, $C^a_b$ must be a permutation matrix 
and the conclusion follows.   
\end{proof}

\begin{definition}[Abstract quantum $D$-module]
$1$. For a given $D$-module $E$, we set 
\[V\defequal E/\left(\sum_{i=1}^r \Qq^i E+\hbar E\right),\quad  
  E_0\defequal E/\sum_{i=1}^r\Qq^iE.
\]
Note that the $D$-submodule $\sum_{i=1}^r \Qq^iE$ does not 
depend on a choice of generalized coordinates of $D$. 
We call $E_0$ the {\bf zero fiber} of $E$. 
$V$ is a module over the commutative algebra $\C[p_1,\dots,p_r]$ 
and $E_0$ is a module over $\C[p_1,\dots,p_r,\hbar]$, 
where the $p_a$'s are the commuting variables corresponding to 
the $\Qp_a$'s. 

$2$. An {\bf abstract quantum $D$-module} is a $D$-module $E$ 
endowed with a base point $e_0$ in $V$ satisfying the following axioms.  
\begin{enumerate}
\item $V$ is a finite dimensional $\C$-vector space. 
\item There exists a splitting  
      $\Phi\colon V\otimes \mathcal{O}^\hbar\rightarrow E$ 
      such that $\Phi$ is an isomorphism of $\mathcal{O}^\hbar$-modules. 
      We call $\Phi$ a {\bf frame} of $E$. 
\item Passing to the quotient by the submodule generated 
      by the action of $q^a$'s, 
      we have an induced isomorphism 
      $\Phi_0\colon V\otimes \C[\hbar]\rightarrow E_0$ 
      from $\Phi$. 
      Then, $\Phi_0$ is an isomorphism of $\C[p_1,\dots,p_r,\hbar]$-modules.
      We call $\Phi_0$ a frame of $E_0$. 
\item The subset $\{e_0,p_1e_0,\dots,p_re_0\}$ of $V$ is linearly independent.
      (Note that this property is independent of a choice of 
       generalized coordinates of $D$.) 
\end{enumerate}  
\end{definition}

It is easy to see that the quantum $D$-module arising from 
quantum cohomology satisfies the above axioms for 
$e_0=1\in V=H^*(M,\C)$ and $\Phi=\id$. 

We postulate only the existence of the frame $\Phi$, 
and do not fix a choice of it. 
When we have two frames $\Phi$ and $\hPhi$, 
they are connected by a {\bf gauge transformation} $Q$ in 
$\Aut_{\mathcal{O}^\hbar}(V\otimes\mathcal{O}^{\hbar})$ such that 
$\hPhi=\Phi\circ Q$. 
The frame of $E_0$ is not fixed either. 
The two frames $\Phi_0$, $\hPhi_0$ of $E_0$ induced by 
$\Phi$, $\hPhi$ are connected also by the gauge transformation 
$Q_0\defequal Q|_{q=0}$.   
$Q_0$ must be an isomorphism of $\C[\hbar,p_1,\dots,p_r]$-modules.   
Moreover, $Q_0|_{\hbar=0}$ must be the identity operator of $V$. 

In this paper, we are interested in the case where 
$V$ is generated by $e_0$ as a $\C[p_1,\dots,p_r]$-module. 
In the case of the quantum cohomology of $M$, 
this condition means that $H^*(M)$ is generated by $H^2(M)$ as a ring.
In this case, a gauge transformation 
$Q_0$ of $E_0$  is determined by the value 
$Q_0(e_0)=e_0+a_1\hbar+a_2\hbar^2+\cdots  (a_i\in V)$ 
which is a generator of $V\otimes\C[\hbar]$ 
as a $\C[\hbar,p_1,\dots,p_r]$-module. 
Therefore, we have a one-to-one correspondence
\[ \text{a frame $\Phi_0$  of }  E_0
   \Longleftrightarrow 
   \text{a generator of }E_0 
   \text{ which projects to }e_0\in V.
\]
In the case of the original quantum $D$-modules, 
there is a natural grading on $E_0$ 
such that $\deg \hbar=2$.  
Therefore, the lift of $e_0$ to $E_0$ can be 
canonically determined by the homogeneity condition 
$\deg\Phi_0(e_0)=0$.

Take a frame $\Phi$ and a set of generalized coordinates 
$(\Qq^1,\dots,\Qq^r,\Qp_1,\dots,\Qp_r)$ of $D$. 
We define a flat connection $\nabla$ 
on the sections of the trivial sheaf 
$V\otimes\mathcal{O}^{\hbar}$ by 
\[ \nabla_a\defequal \nabla_{\parti{a}}
\defequal \frac{1}{\hbar}\text{ the action of $\Qp_a$} 
\]
through the identification 
$\Phi\colon V\otimes\mathcal{O}^\hbar\cong E$. 
This is considered to be a connection 
over the infinitesimal neighborhood of 0 in $B$. 
More specifically, 
define connection matrices 
$\Omega_a$ in $\End(V)\otimes\mathcal{O}^{\hbar}$  
for $a=1,\dots,r$ by  
\[\Qp_a\cdot\Phi(v)=\Phi(\Omega_a(v)), \quad v\in V.
\]
The flat connection $\nabla$ is defined by  
\[\nabla\colon V\otimes\mathcal{O}^{\hbar}\longrightarrow 
 \frac{1}{\hbar}V\otimes\sum_{a=1}^r\mathcal{O}^{\hbar}
 \frac{dq^a}{q^a}, \quad 
  \nabla=d+\frac{1}{\hbar}\sum_{a=1}^r \Omega_a \frac{dq^a}{q^a}.
\]
In particular, 
$\nabla_{a}(v)\defequal \frac{1}{\hbar}\Omega_a(v)$ for 
$v\in V$. 
The connection matrix $\Omega_a(q=0,\hbar)$ at the origin 
is $\hbar$-independent 
and is equal to the action of $p_a$ on $V$.  
This follows from the third axiom. 
By the relation $[\Qp_a,\Qp_b]=0$, the connection $\nabla$ is flat.  
It corresponds to 
the dual Givental connection $\nabla^{\hbar}$ divided by $\hbar$ 
in the original quantum $D$-module case. 
We describe the transformation law of connection matrices. 
Let $\Phi$ and $\hPhi=\Phi\circ Q$ be frames. 
The corresponding connection matrices 
$\Omega_a$ and $\hOmega_a$ are related by  
\begin{equation}
\label{eq:connectionchangeviagaugetransformation}
\hOmega_{a}=Q^{-1}\Omega_aQ+Q^{-1}\hbar\parti{a}Q.
\end{equation}
The connection matrices are transformed also
under the coordinate transformation (\ref{eq:coordinatetransformation})
as follows:   
\begin{equation}
\label{eq:connectionchangeviamirrortransformation} 
\Omega_{\hat{a}}=\sum_{b=1}^r\lparfrac{\log q^b}{\hq^a}\Omega_b+
                    \lparfrac{F}{\hq^a}.
\end{equation}
Note that connection matrices are transformed 
in two different ways.

\subsection{Canonical Frame}

A frame $\Phi$ of $E$ is called {\bf canonical} 
if its connection matrices $\Omega_a$ are $\hbar$-independent 
i.e. $\Omega_a\in \End(V)\otimes\mathcal{O}$. 
The notion of canonical frame does not depend on 
a choice of generalized coordinates of $D$. 
In this subsection,            
we show that given a frame $\Phi_0$ of $E_0$, 
there exists a unique canonical frame $\Phi$ 
which induces $\Phi_0$.  

For the original quantum $D$-module, the canonical frame is 
already chosen by the canonical isomorphism 
$E\cong H^*(M)\otimes \C[\hbar][\![q]\!]$. 
Connection matrices are identified with quantum products 
by two dimensional classes, 
therefore they do not depend on $\hbar$. 
The existence and uniqueness of the canonical frame shows that 
the abstract $D$-module structure 
reconstructs the quantum multiplication table. 
(In this case, the frame of $E_0$ is determined by homogeneity.)  

We introduce flat connections $\nabla^0$ and $\nabla^1$ 
on the endomorphism bundle 
$\End(V)\otimes \mathcal{O}^{\hbar}$. 
Given a frame of $E$ and generalized coordinates of $D$, 
they are defined by
\begin{equation}
\label{eq:connection01}
  \nabla^0_aT\defequal \parti{a}T+\frac{1}{\hbar}\Omega_aT,
  \quad 
  \nabla^1_aT\defequal \parti{a}T+\frac{1}{\hbar}(\Omega_aT-Tp_a),  
\end{equation} 
for $T\in \End(V)\otimes\mathcal{O}^{\hbar}$, 
where $p_a$ is considered to be an element of $\End(V)$.  
It is easy to check that they are flat. 
The connection $\nabla^0$ defines a fundamental solution $L$ 
for $\nabla$ which is an analogue of the $L$ in 
(\ref{eq:non-abstractL}). 
A parallel section $S$ for the connection $\nabla^1$ 
is an analogue of the $S$ in (\ref{eq:non-abstractS}). 
In fact, the differential equation in 
Lemma \ref{lem:non-abstractdiffeqforS} is identical with that 
for $\nabla^1$. 
First we construct a solution $S$ for $\nabla^1$, 
and next we define a fundamental solution for $\nabla$ 
by $L\defequal S\circ e^{-p\log q/\hbar}$, 
where $p\log q=\sum_{a=1}^rp_a\log q^a$. 
Note that any solution $\tilde{L}$ for $\nabla^0$ regular at $q=0$ 
satisfies $p_a\tilde{L}(q=0)=0$ and hence is not invertible. 
Therefore, an invertible matrix solution $L$
for $\nabla^0$ cannot be regular at $q=0$.

\begin{proposition} 
\label{prop:existenceofS}
There exists a unique $\nabla^1$-flat section $S(q,\hbar)$ 
in $\End(V)\otimes\C[\hbar,\hbar^{-1}]\!][\![q^1,\dots,q^r]\!]$
which satisfies the initial condition $S(0,\hbar)=\id$.   
The solution $S(q,\hbar)$ depends on a choice of 
a frame of $E$ and generalized coordinates of $D$. 
\end{proposition} 
\begin{proof}
We solve the second equation of (\ref{eq:connection01})
with the initial condition $S(0,\hbar)=\id$.  
Suppose by induction that we have a solution 
$S(k)\defequal S(q^1,\dots,q^k,0,\dots,0,\hbar)$ which 
satisfies the second equation of (\ref{eq:connection01}) for $a\le k$. 
We assume that 
\[\Omega_a(q^1,\dots,q^k,0,\dots,0,\hbar)S(k)
=S(k)p_a,
\quad  \forall a>k.
\] 
These assumptions are satisfied when $k=0$. 
We expand $S$ and 
$\Omega_{k+1}$ as follows. 
\begin{gather*} 
 S(k+1)=S(k) \sum_{n=0}^\infty T_n(q^1,\dots,q^k,\hbar)(q^{k+1})^n, \\
 \Omega_{k+1}(q^1,\dots,q^{k+1},0,\dots,0,\hbar)
       =\sum_{n=0}^\infty\Omega_{k+1,n}(q^1,\dots,q^{k},\hbar)(q^{k+1})^n, 
\end{gather*} 
where $T_0=\id$. 
We substitute $S(k+1)$ for the second equation of 
(\ref{eq:connection01}) in the case $a=k+1$ and expand in 
$q^{k+1}$. 
Extracting the coefficient of $(q^{k+1})^n$, 
we obtain for $n\ge 1$, 
\[nS(k)T_n+\frac{1}{\hbar}\left\{\sum_{i=0}^{n-1}\Omega_{k+1,n-i}S(k)T_i
                           +\Omega_{k+1,0}S(k)T_n 
                           -S(k)T_np_{k+1}\right\}=0. \]
By the assumption $\Omega_{k+1,0}S(k)=S(k)p_{k+1}$, 
we have for $n\ge 1$, 
\[\left(n+\frac{1}{\hbar}\ad(p_{k+1})\right)T_n
       =-\frac{1}{\hbar}\sum_{i=0}^{n-1}S(k)^{-1}\Omega_{k+1,n-i}S(k)T_i. 
\]
If $T_1,\dots,T_{n-1}$ are given, 
we can solve for $T_n$ in the above equation. 
Note that $T_n$ is contained in the ring 
$\End(V)\otimes\C[\hbar,\hbar^{-1}]\!][\![q^1,\dots,q^k]\!]$. 
Thus we obtain $S(k+1)$.

Next, we verify that the second equation 
of (\ref{eq:connection01}) 
holds for $S=S(k+1)$ and $a=1,\dots,k$. 
For $1\le a\le k$, we have 
\[\nabla^1_{k+1}\nabla^1_{a}S(k+1)
 =\nabla^1_{a}\nabla^1_{k+1}S(k+1)=0, \quad
 \nabla^1_aS(k+1)|_{q^{k+1}=0}=0. 
\]
Hence $T(q^1,\dots,q^{k+1},\hbar)\defequal\nabla^1_{a}S(k+1)$ 
satisfies the differential equation $\nabla^1_{k+1}T=0$ and 
the initial condition $T|_{q^{k+1}=0}=0$. 
This equation turns out to have a unique solution and 
we have $\nabla^1_aS(k+1)=T\equiv 0$.    

Finally we check that for $a>k+1$, the relation 
\[\Omega_a(q^1,\dots,q^{k+1},0,\dots,0,\hbar)S(k+1)-S(k+1)p_a=0\]
holds.  
Because the left hand side 
is equal to $\nabla^1_aS(k+1)|_{q^{k+2}=\dots=q^{r}=0}$,  
it turns out to be zero by the same argument.  
This completes the induction step, and we obtain
$S(q,\hbar)$.  
\end{proof} 

\begin{remark}
\label{rem:regularityathbarinfty}
In the case where all $\Omega_a$ are $\hbar$-independent, 
the section $S(q,\hbar)$ obtained in the previous proposition 
is regular at $\hbar=\infty$. 
In other words, $S(q,\hbar)\in \C[\![h^{-1}]\!][\![q^1,\dots,q^r]\!]$.  
\end{remark}

We describe the dependency of the solution $S(q,\hbar)$
in Proposition \ref{prop:existenceofS} 
on a frame and generalized coordinates.  
Denote by $S_{\Phi,(\Qq,\Qp)}$ the solution 
corresponding to a frame $\Phi$ of $E$ 
and generalized coordinates $(\Qq,\Qp)$ of $D$. 
\begin{lemma}
\label{lem:transformationofS}
Let $\Phi$ and $\hPhi$ be frames of $E$ and 
$Q(q,\hbar)$ be the gauge transformation such that $\hPhi=\Phi\circ Q$. 
Define $Q_0\defequal Q(0,\hbar)$.  
Let $(\Qq,\Qp)$ and $(\hQq,\hQp)$ be two sets of 
generalized coordinates of $D$.  
The two solutions $S_{\Phi,(\Qq,\Qp)}$ and $S_{\hPhi,(\hQq,\hQp)}$ 
are related by 
\[S_{\hPhi,(\hQq,\hQp)}=Q^{-1}S_{\Phi,(\Qq,\Qp)}Q_0
  e^{\sum_{a=1}^rp_a(\log \hq^a-\log q^a-c^a)/\hbar-F/\hbar},\]
where $F$ is an element of $\mathfrak{m}$ associated with the
coordinate transformation which appeared in the equation 
$(\ref{eq:coordinatetransformation})$ and 
$c^a\in \C$ is chosen so that 
$\log \hq^a-\log q^a-c^a$ is in $\mathfrak{m}$. 
(Note that $\hq^a$ is of the form 
$\hq^a=q^a \exp(\delta^a),\; \delta^a\in \mathcal{O}$ 
after renumbering the indices.
 See Lemma \ref{lem:integrabilitycondforlog-Jacobi}. )
\end{lemma}
\begin{proof}
The initial condition $S_{\hPhi,(\hQq,\hQp)}(\hq=0,\hbar)=\id$ 
is easy to check. 
Using the fact that $Q_0$ commutes with $p_a$ and
the equations  
(\ref{eq:connectionchangeviagaugetransformation}), 
(\ref{eq:connectionchangeviamirrortransformation}),  
we can check that the left hand side satisfies the 
differential equation. 
\end{proof}

We prove the main theorem of this section,  
the existence and uniqueness of the canonical frame. 
The theorem is stated as follows. 
\begin{theorem} 
\label{thm:canonicalframe}
For a given frame $\Phi_0$ of $E_0$, 
there exists a unique canonical frame $\hPhi$ of $E$ 
which induces $\Phi_0$.  
In other words, for any frame $\Phi$ which induces $\Phi_0$, 
there exists a unique gauge transformation $Q(q,\hbar)$ such that 
$Q(0,\hbar)=\id$ and $\hPhi\defequal\Phi\circ Q(q,\hbar)$ 
is canonical. 
\end{theorem} 
\begin{proof}  
Fix a set of generalized coordinates $(\Qq,\Qp)$ of $D$ and 
choose a frame $\Phi$ of $E$.  
We use the flat section $S(q,\hbar)$ of $\nabla^1$ 
obtained in Proposition \ref{prop:existenceofS}. 
First we perform Birkhoff factorization for $S(q,\hbar)$. 
We prove that $S$ can be uniquely factorized as 
$S=S_{+}S_{-}$, where $S_{+}\in \End(V)\otimes\C[\hbar][\![q]\!]$ 
and $S_{-}\in \End(V)\otimes\C[\![\hbar^{-1}]\!][\![q]\!]$ 
and $S_{-}(\hbar=\infty)=\id$. 
We expand $S$, $S_{+}$ and $S_{-}$ 
in formal power series of $q^a$ as follows. 
\[S=\id+\sum_{d\neq 0} S_d(\hbar,\hbar^{-1})q^d,  \quad
  S_{+}=A_0(\hbar)+\sum_{d\neq 0} A_d(\hbar)q^d, \quad
  S_{-}=\id+\sum_{d\neq 0}B_d(\hbar^{-1})\frac{1}{\hbar}q^d, 
\]
where $d\in(\Z_{\ge 0})^r$, $A_d\in \C[\hbar]$ 
and $B_d\in \C[\![\hbar^{-1}]\!]$. 
We determine $A_d$,  $B_d$ inductively. 
It is easy to see that $A_0=\id$. 
We define $d_1\le d_2$ for $d_i=(d_{i,1},\dots,d_{i,r})\in (\Z_{\ge 0})^r$ 
if $d_{1,j}\le d_{2,j}$ holds for all $j$. 
Assume that we obtain $A_d$ and $B_d$ for all $d<d_0$. 
Then we have 
\[S_{d_0}(\hbar,\hbar^{-1})=B_{d_0}(\hbar^{-1})\frac{1}{\hbar}+A_{d_0}(\hbar)+
\sum_{\substack{d_1+d_2=d_0 \\ d_1\neq 0,d_2\neq 0}} 
 A_{d_1}(\hbar)B_{d_2}(\hbar^{-1})\frac{1}{\hbar}.
\] 
This determines $A_{d_0}$ and $B_{d_0}$ uniquely. 
Note that $S_{+}(0,\hbar)=S_{-}(0,\hbar)=\id$. 
By the differential equation $\nabla^1S(q,\hbar)=0$, we have
\[S^{-1}\parti{a}S+S^{-1}\frac{\Omega_a}{\hbar}S=\frac{p_a}{\hbar}. 
\] 
Substituting $S_{+}S_{-}$ for $S$, we obtain 
\[S_{+}^{-1}\parti{a}S_{+}+S_{+}^{-1}\frac{\Omega_a}{\hbar}S_{+}=
            S_{-}\frac{p_a}{\hbar}S_{-}-(\parti{a}S_{-})S_{-}^{-1}.\]  
The $\hbar$-expansion of the left hand side is 
$A^{-1}\Omega|_{\hbar=0}A/\hbar+O(1)$ for $A=S_{+}(\hbar=0)$ and 
the right hand side is of order $O(1/\hbar)$. 
Therefore, we obtain 
\[S_{+}^{-1}\parti{a}S_{+}+S_{+}^{-1}\frac{\Omega_a}{\hbar}S_{+}
          =\frac{\hOmega_a}{\hbar}, \qquad 
          \text{for some $\hbar$-independent matrix } 
          \hOmega_a. 
\]
The above equation shows that 
$S_{+}$ gives the desired gauge transformation $Q$.  

Next, we show that the gauge transformation $Q$ is unique. 
Assume that for another gauge transformation $Q'$ satisfying 
$Q'(q=0)=\id$, we have a new frame $\Phi'=\Phi\circ Q'$ 
such that its connection matrices $\Omega_a'$ are $\hbar$-independent. 
For the frame $\Phi'$, take the unique solution $S'(q,\hbar)$ 
for the following system of differential equations. 
\[{\nabla^1_a}'S'(q,\hbar)=\parti{a}S'(q,\hbar)+
  \frac{1}{\hbar}(\Omega_a'S'-S'p_a)=0, \quad
   S'(0,\hbar)=\id
\]  
By Remark \ref{rem:regularityathbarinfty}, 
$S'$ is regular at $\hbar=\infty$. 
Set $\tilde{S}(q,\hbar)=Q'(q,\hbar)S'(q,\hbar)$. 
Then, it is easy to show that 
$\nabla^1\tilde{S}=0$ and $\tilde{S}(0,\hbar)=\id$. 
By the uniqueness of the solution $S$, 
we have $S=\tilde{S}$. 
Also by the uniqueness of Birkhoff factorization, 
we have $Q=S_{+}=Q'$.  
\end{proof}

\begin{remark}
We owe the idea of using the Birkhoff factorization to Martin Guest 
\cite{guest-Dmodule} in this $D$-module context. 
He used this factorization not for $S$, 
but for any fundamental solution matrix for $\nabla$. 
The Birkhoff factorization also appears in the works of 
Barannikov \cite{barannikov} and Coates and Givental \cite{coates-givental}. 
\end{remark} 
\begin{remark}
The computation of the canonical frame and its connection can be carried out 
in an algorithmic way. 
We can compute the coefficient of each $q^d$ step by step, 
but we cannot obtain a closed formula by this method.    
\end{remark}

As an application of Theorem \ref{thm:canonicalframe},   
we show that the $J$-function in Proposition \ref{prop:explicitJ} 
determines all genus zero Gromov-Witten invariants of $M$ 
if $H^*(M)$ is generated by $H^2(M)$. 
By Theorem \ref{thm:Jisgen}, we can reconstruct 
the quantum $D$-module from the $J$-function. 
In this case, the image of $1\in D$ in $D/I$ gives 
the base point $\Phi(e_0)$ of the quantum $D$-module, 
therefore, the frame of the zero fiber $E_0$ is already given.   
We remarked that the quantum $D$-module arising from the quantum cohomology 
has connection matrices independent of $\hbar$. 
Therefore, taking the canonical frame, 
we have a table of quantum multiplication by 
two dimensional classes.  
By the reconstruction theorem of Kontsevich and Manin, 
we can determine all genus zero Gromov-Witten invariants
\cite{kontsevich-manin}. 

\begin{corollary}
\label{cor:JfunctiondeterminesGW}
Let $(M,\mV)$ be a supermanifold which satisfies 
the assumptions in Theorem \ref{thm:Jisgen}. 
The $J$-function defined by the small quantum cohomology of $(M,\mV)$ 
reconstructs not only the isomorphism type of the quantum cohomology ring 
but also the quantum multiplication table. 
Moreover, when $\mV=0$, 
it determines all genus zero Gromov-Witten invariants of $M$.  
\end{corollary}

\subsection{The Graded Case}

The original quantum $D$-module has a natural grading. 
We introduce the notion of grading for abstract quantum $D$-modules. 
In the graded case, the connections $\nabla^0$, $\nabla^1$ 
can be extended over $B\times \C$, 
where the coordinate of $\C$ is $\hbar$.  
Further, we can take coordinates compatible  
with canonical frame if the grading 
satisfies the nef condition 
(Assumption \ref{ass:nef}).  

In the graded case, 
we assume that 
the variables have the following degrees.  
\[\deg \Qp_a=\deg p_a=\deg\hbar=2,\quad \deg \Qq^a=\deg q^a=
  \text{ even integer}.
\] 
We only consider coordinate transformations 
which preserve the homogeneity of the degree. 
We define the {\bf Euler vector field} over $B\times \C$ 
to be 
\[\Eul=2\hbar\parfrac{}{\hbar}+\sum_{a=1}^r(\deg q^a)\parti{a}.\] 
We define an operator $\adE:D\rightarrow D$ satisfying 
the Leibniz rule as follows:  
\begin{gather*}
  \adE(\Qp_a)=2\Qp_a,\quad \adE(\hbar)=[\Eul,\hbar]=2 \hbar, 
  \quad \adE(\Qq^a)=[\Eul,\Qq^a]=(\deg q^a)\Qq^a, \\
  \adE(x\cdot y)=\adE(x)\cdot y+x\cdot\adE(y) \quad 
  \text{ for $x$, $y$ $\in D$.} 
\end{gather*} 

The operator $\adE$ acts on the abstract quantum $D$-module $E$ 
and also satisfies the Leibniz rule: 
\begin{equation}
\label{eq:actionofadE}
\adE(x\cdot e)=\adE(x)\cdot e+x\cdot\adE(e),  
\end{equation}
for $x\in D$ and $e\in E$. 
The action of $\adE$ on $E$ induces an action on $V$. 
If $\adE$ is semisimple on $V$, its eigenvalues 
are the degrees of the eigenvectors. 

We assume in the graded case that there exists a frame 
$\Phi$ of $E$ which commutes with the action of $\adE$, i.e. 
$\Phi\circ\adE=\adE\circ\Phi$, 
where the action of $\adE$ is extended on $V\otimes\mathcal{O}^\hbar$ 
in such a way that it satisfies equation (\ref{eq:actionofadE})
for $x\in \mathcal{O}^\hbar$ and $e\in V$. 
We only consider frames commuting with $\adE$. 
Therefore, the gauge transformation $Q(q,\hbar)$ and 
the frame $\Phi_0$ of $E_0$ also commute with $\adE$.  
Let $\mu\in \End(V)$ be the action of $\adE$ on $V$ and 
put $\rho\defequal \sum_{a=1}^r(\deg q^a)p_a\in \End(V)$. 
For a given frame $\Phi$ and generalized coordinates $(\Qq,\Qp)$, 
we extend the connections $\nabla^0$, $\nabla^1$ 
over $B\times \C$
by the following formula.  
\begin{eqnarray*}
\nabla^0_{\Eul}T &\defequal & [\adE, T]+T\frac{\rho}{\hbar}
=\Eul T+[\mu,T]+T\frac{\rho}{\hbar}, \\
\nabla^1_{\Eul}T &\defequal & [\adE, T]=\Eul T+[\mu,T], 
\end{eqnarray*}
for $T\in \End(V)\otimes \mathcal{O}^\hbar$. 
The extended connections $\nabla^0$, $\nabla^1$ 
(which we denote by the same symbols) are also flat. 
In the case of quantum cohomology, 
Lemma \ref{lem:non-abstractdiffeqforS} shows that 
$S$ in equation (\ref{eq:non-abstractS}) 
satisfies $\nabla^1_{\Eul}S=0$. 
If $S$ satisfies $\nabla^1_{\Eul}S=0$, 
$L=S\circ e^{-p\log q/\hbar}$ satisfies 
$\nabla^0_{\Eul}L=0$. 
We show that there exists a unique solution $S$ for 
the extended connection $\nabla^1$ with the initial condition 
$S(q=0,\hbar)=\id$. 
In other words,  
\begin{proposition}
The solution $S$ constructed in 
Proposition \ref{prop:existenceofS} 
also satisfies $\nabla^1_{\Eul}S=0$ in the graded case. 
Therefore, $S$ automatically preserves the degree. 
Moreover, $S$ is an element of 
$\End(V)\otimes\C[\hbar,\hbar^{-1}][\![q^1,\dots,q^r]\!]$. 
\end{proposition}
\begin{proof}
Let $S(q,\hbar)$ be the solution in Proposition \ref{prop:existenceofS}. 
By the flatness of $\nabla^1$, we have 
\[\nabla^1_a\nabla^1_{\Eul}S=\nabla^1_{\Eul}\nabla^1_aS=0.\]
Moreover, $\nabla^1_{\Eul}S$ satisfies the initial condition 
$(\nabla^1_{\Eul}S)(q=0,\hbar)=0$. 
From this it follows that $\nabla^1_{\Eul}S=0$. 
The last statement follows from $[\adE,S]=0$,  
$\deg \hbar=2$ and $\dim_{\C}V<\infty$. 
\end{proof}

By this proposition, Theorem \ref{thm:canonicalframe} 
also holds in the graded case. 
\begin{theorem}
For a given frame $\Phi_0$ of $E_0$ commuting with $\adE$, 
there exists a unique canonical frame $\hPhi$ which induces $\Phi_0$ 
and commutes with $\adE$. 
In other words, for any frame $\Phi$ which induces 
$\Phi_0$ and commutes with $\adE$, there exists 
a unique degree-preserving gauge transformation $Q(q,\hbar)$ 
such that $Q(0,\hbar)=\id$ and $\hPhi\defequal \Phi\circ Q(q,\hbar)$
is canonical.  
\end{theorem} 
\begin{proof}
Let $S$ be the solution for $\nabla^1$ obtained in the
previous proposition. 
It suffices to check that the positive part $S_+$ of the 
Birkhoff factorization $S=S_+S_-$ preserves the degree. 
Because $\nabla^1_{\Eul}S=[\adE,S]=0$, we have 
\[[\adE,S_-]S_-^{-1}=-S_+^{-1}[\adE,S_+].
\]
Because the left hand side is of order $O(1/\hbar)$, 
both sides are zero. 
\end{proof}

In the case of original quantum $D$-modules, 
we have the relation $\Qp_a\cdot\Phi(e_0)=\Phi(p_ae_0)$ 
for $e_0=1\in V=H^*(M,\C)$. 
In the case of an abstract quantum $D$-module, 
we call generalized coordinates $(\Qq^1,\dots,\Qq^r,\Qp_1,\dots,\Qp_r)$ 
{\bf compatible with a frame $\Phi$} 
if they satisfy the same relation 
$\Qp_a\cdot\Phi(e_0)=\Phi(p_a e_0)$.  
In general, we can not expect that 
compatible generalized coordinates of $D$ exist. 
However, if they exist, 
they are uniquely determined up to the ordering and 
multiplication by constants $\Qq^a\mapsto c^a\Qq^a$. 
This follows from the coordinate transformation 
(\ref{eq:coordinatetransformation}),  
the fact that $\{e_0,p_1e_0,\dots,p_re_0\}$ is linearly independent, 
and Lemma \ref{lem:integrabilitycondforlog-Jacobi}. 
Next, we give a sufficient condition 
for the existence of generalized coordinates compatible 
with a canonical frame. 
For that, we assume 
\begin{assumption}
\label{ass:nef}
$1.$ $\mu\in \End(V)$ is semisimple and 
$V$ decomposes into the direct sum 
\[ V=V_0\oplus V_2 \oplus V_4 \oplus \dots \oplus V_{2n}, \]
where  $V_{2k}$ is the eigenspace of $\mu$ with eigenvalue $2k$.  

$2.$ $V_0=\C e_0$ and $V_2=\bigoplus_{a=1}^r\C (p_ae_0)$. 

$3.$ The grading satisfies the nef condition: $\deg q^a\ge 0$ 
    for $a=1,\dots,r$. 
\end{assumption}  
The quantum $D$-module arising from the quantum cohomology 
of the superspace $(M,\mV)$ satisfies this assumption 
if its first Chern class $c_1(M/\mV)$ is nef. 

\begin{theorem}
\label{thm:mirrortransformation}
Under Assumption \ref{ass:nef}, there exist generalized coordinates  
$(\Qq,\Qp)$ compatible with a canonical frame $\Phi$ of $E$. 
The set of compatible generalized coordinates $(\Qq,\Qp)$ 
is unique up to the ordering and multiplication by constants  
$\Qq^a\mapsto c^a\Qq^a$.
\end{theorem} 
\begin{proof}
Let $\Omega_a$ be a connection matrix associated with a 
canonical frame $\Phi$. 
By Assumption \ref{ass:nef} and 
$\deg \Omega_a=2$ (i.e. $[\adE,\Omega_a]=2\Omega_a$), 
one can write 
\begin{equation}
\label{eq:mirrortransmethod1} 
\Omega_a(e_0)=-F_a(q)e_0
  +\sum_{b=1}^r G^b_a(q)(p_be_0), \quad 
  F_a, G^b_a\in \mathcal{O}, 
\end{equation} 
where $F_a(q=0)=0$ and $G^b_a(q=0)=\delta^b_a$. 
By the flatness of $\nabla=d+(1/\hbar)\sum_{a=1}^r\Omega_a dq^a/q^a$, 
we have 
\[\parti{a}\Omega_b-\parti{b}\Omega_a
 +\frac{1}{\hbar}[\Omega_a,\Omega_b]=0.
\]
Because $\Omega_a$ is $\hbar$-independent, we have 
\[\parti{a}\Omega_b-\parti{b}\Omega_a=0,\quad 
  [\Omega_a,\Omega_b]=0.
\]
Therefore, by Lemma \ref{lem:integrabilitycond} and  
Lemma \ref{lem:integrabilitycondforlog-Jacobi},  
there exists a function $F\in\mathfrak{m}$ and 
a new set of coordinates $\{\hq^1,\dots,\hq^r\}$ 
such that 
\begin{equation}
\label{eq:mirrortransmethod2} 
F_a(q)=\parti{a}F(q),\quad G^b_a(q)=\lparfrac{\log\hq^b}{q^a}. 
\end{equation} 
Define new generalized coordinates $\hQp_a$ by 
\[\hQp_a=\sum_{b=1}^r \lparfrac{\log q^b}{\hq^a}\Qp_b 
         +\parti{a}F(q). 
\]
Then equation (\ref{eq:connectionchangeviamirrortransformation}) 
shows that $\Omega_{\hat{a}}(e_0)=p_ae_0$, 
therefore $(\hQq,\hQp)$ are compatible generalized coordinates. 
\end{proof} 

The coordinate transformation appearing in the above theorem 
is called a {\bf mirror transformation}. 
By choosing generalized coordinates compatible with a canonical frame, 
we can canonically define an affine structure on the base $B$. 
(The coordinates $\log q^a$ define a flat structure.)
The equivariant Floer cohomology introduced later 
is not a priori given the `correct' (affine) coordinates. 
We can find an affine structure by taking the generalized
coordinates compatible with the canonical frame. 

The indeterminacy of the constant multiple $q^a\mapsto c^aq^a$ 
can be normalized by choosing an element $v_0$ of 
$T_0B^*\defequal
\{\sum_{a=1}^r c^a(\partial/\partial q^a)_0
     |\; c^a\neq 0\;(\forall a)\}$.
Then we normalize $q^a$ so that the equality 
$v_0=(\partial/\partial q^1)_0+\dots+(\partial/\partial q^r)_0$ 
holds. 
We call $v_0\in T_0B^*$ a {\bf normalization vector}. 
In this sense, the generalized coordinates in the equivariant Floer theory 
are a priori correctly normalized 
even though they are not affine. 

\subsection{$J$-function and Pairing}

In this subsection, we define the $J$-function and 
the pairing of the abstract quantum $D$-module. 
In what follows, we fix a choice of a frame $\Phi_0$ of 
the zero fiber $E_0$ and a normalization vector $v_0\in T_0B^*$. 
The $J$ function depends on a choice of a frame, but 
the pairing does not depend on it. 

Let $\Phi$ be a frame of $E$ and 
$(\Qq,\Qp)$ be a set of generalized coordinates of $D$. 
The flat section $L_{\Phi,(\Qq,\Qp)}$ 
of $\nabla^0$ can be obtained from 
the solution $S_{\Phi,(\Qq,\Qp)}$ in Proposition \ref{prop:existenceofS} 
by the formula
\[L_{\Phi,(\Qq,\Qp)}\defequal 
  S_{\Phi,(\Qq,\Qp)}\circ e^{-p\log q/\hbar}, 
  \quad p\log q=\sum_{a=1}^rp_a\log q^a. 
\]
The function $L_{\Phi,(\Qq,\Qp)}$ is 
a multi-valued matrix function on $B$. 
By Lemma \ref{lem:transformationofS},  
we have the following transformation law: 
\[L_{\hPhi,(\hQq,\hQp)}=e^{-F/\hbar}Q^{-1}L_{\Phi,(\Qq,\Qp)}. 
\]
Note that $Q_0=\id$ and $\log q^a-\log \hq^a\in \mathfrak{m}$ 
because we fixed a choice of a frame of the zero fiber 
and a normalization vector. 
 
Next, we define the $J$-function as follows.  
\[J_{\Phi,(\Qq,\Qp)}\defequal L_{\Phi,(\Qq,\Qp)}^{-1}(e_0)
=e^{p\log q/\hbar}S_{\Phi,(\Qq,\Qp)}^{-1}(e_0).\]  
The following theorem is the analogue of Theorem \ref{thm:Jisgen}. 
\begin{theorem} 
\label{thm:abstractreconstructionbyJ}
Assume that $V$ is generated by $e_0$ as a $\C[p_1,\dots,p_r]$-module. 
Then $E$ is generated by a single element $\Phi(e_0)$ as a $D$-module. 
Moreover, we have a $D$-module isomorphism $D/I\cong E$ which sends 
$1$ to $\Phi(e_0)$, where 
$I$ is the left ideal consisting of elements $f(\Qq,\Qp,\hbar)$ in 
$D$ satisfying 
\[f(q,\hbar\parti{{}},\hbar)\cdot J_{\Phi,(\Qq,\Qp)}=0.
\]  
\end{theorem} 
The proof is similar to that of Theorem \ref{thm:Jisgen}. 

When we have a canonical frame and compatible generalized coordinates, 
we can describe the asymptotic behavior of the $J$-function. 
By Remark \ref{rem:regularityathbarinfty}, 
the solution $S_{\rm can}=S_{\Phi_{\rm can},(\Qq,\Qp)}$ 
can be expanded in $\hbar^{-1}$ and one can write 
$S_{\rm can}^{-1}=S_0+\hbar^{-1}S_1+\hbar^{-2}S_2+\cdots$. 
From $\nabla^1S_{\rm can}=0$, we have 
\[\parti{a}S_0=0, \qquad 
\parti{a}S_1+p_aS_0-S_0\Omega_a^{\rm can}=0.\]
Together with the initial condition $S_{\rm can}|_{q=0}=\id$ 
and the compatibility $\Omega_a^{\rm can}(e_0)=p_ae_0$,   
these show that $S_0=\id$ and $S_1(e_0)=0$. 
Thus we obtain the following. 
\begin{equation}
\label{eq:propertyofcanonicalJ}
J_{\rm can}=J_{\Phi_{\rm can}}=e^{p\log q/\hbar}(e_0+o(\hbar^{-1})).
\end{equation} 

Under Assumption \ref{ass:nef}, 
a degree-preserving gauge transformation $Q$ satisfies 
$Q(e_0)=f(q)e_0$ for some degree zero element $f\in \mathcal{O}$. 
Therefore, the $J$-function is transformed as 
\begin{equation}
\label{eq:transformationofJ}
J_{\hPhi,(\hQq,\hQp)}=fe^{F/\hbar}J_{\Phi,(\Qq,\Qp)},
\end{equation}
This transformation law seems more natural when we 
take the viewpoint of Remark \ref{rem:senseoft_0} 
and define $\tilde{J}_{\Phi,(\Qq,\Qp)}=e^{t_0/\hbar}J_{\Phi,(\Qq,\Qp)}$. 
Then the term $e^{F/\hbar}$ is absorbed into 
the transformation of $t_0$ and $\tilde{J}_{\Phi,(\Qq,\Qp)}$ is 
transformed as $\tilde{J}_{\hPhi,(\hQq,\hQp)}=f\tilde{J}_{\Phi,(\Qq,\Qp)}$. 
In this case,  
the asymptotic expansion (\ref{eq:propertyofcanonicalJ})
characterizes the canonical $J$-function $J_{\rm can}$ 
among all $J$-functions 
and we can calculate the mirror transformation 
from any $J$-function. 

\begin{proposition}
\label{prop:giventalmirrortrans} 
Let $E$ be a graded abstract quantum $D$-module 
and $J_{\Phi,(\Qq,\Qp)}$ be a $J$-function. 
Under Assumption \ref{ass:nef}, 
we have a unique coordinate transformation of the form, 
\[ \log q^a=\log \hq^a+ \delta^a(\hq), \qquad 
\delta^a(\hq)\in\mathfrak{m}, \]
and $F,f\in \mathcal{O}$, $\deg F=2$, $\deg f=0$, $f(q=0)=1$
satisfying 
\[fe^{F/\hbar}J_{\Phi,(\Qq,\Qp)}=e^{\sum_{a=1}^r p_a\log \hq^a/\hbar}
   (e_0+o(\hbar^{-1})).
\] 
The left hand side gives the canonical $J$-function.   
\end{proposition} 
This proposition corresponds to Givental's mirror theorem concerning   
the relation between two cohomology-valued functions $I$ and $J$. 
In our terminology, $I$ is a $J$-function of a non-canonical 
frame and $J$ is that for the canonical frame. 
Givental showed that a coordinate change of the $I$-function 
gives the $J$-function. 
This proposition is more general than the Givental's original situation 
because this is formulated in terms of abstract quantum $D$-modules.  
We omit the proof because 
it is similar to the proof appearing in  
\cite{cox-katz, givental-mirrorthm-toric}.

Next, we define a pairing for an abstract quantum $D$-module $E$. 
We assume that $V$ has a symmetric pairing 
$\langle \cdot,\cdot\rangle:V\times V\rightarrow \C$ 
satisfying the relation 
\[ \pair{v}{p_a w} =\pair{p_a v}{w}.
\] 
We would like to extend this pairing of $V$ on $E$. 
First we define a $\ovD$-module $\ovE$. 
Let $\ovD$ be a Heisenberg algebra whose commutation 
relations are opposite to $D$: 
\[\ovD=\C[\hbar][\![\Qq^1,\dots,\Qq^r]\!][\Qp_1,\dots,\Qp_r],
\qquad [\Qp_a,\Qq^b]=-\hbar\delta_a^bq^b, 
\quad [\Qp_a,\Qp_b]=[\Qq^a,\Qq^b]=0.\]  
$\ovE$ is a $\ovD$ module which has a bar isomorphism 
$\barop:E\overset{\cong}{\rightarrow} \ovE$ such that 
$\hbar\cdot\ov{x}=\ov{-\hbar\cdot x}$, 
$\Qp_a\cdot\ov{x}=\ov{\Qp_a\cdot x}$, 
$\Qq^a\cdot\ov{x}=\ov{\Qq^a\cdot x}$. 
Briefly, $\ovE$ is $E$ with $\hbar$ action of opposite sign. 
A frame $\Phi$ of $E$ induces a frame $\ovPhi$ of $\ovE$: 
$\ovPhi=\barop\circ\Phi\circ\barop\colon 
V\otimes \mathcal{O}[\hbar]\rightarrow \ovE$, 
where $\barop$ is defined on $V\otimes\mathcal{O}[\hbar]$ 
by $\ov{\hbar}=-\hbar$, $\ov{q^a}=q^a$, $\ov{v}=v$ for $v\in V$.  
The connection $\ov{\nabla}$ associated with the frame 
$\ovPhi$ is given by 
\[\ov{\nabla}=
  d-\frac{1}{\hbar}\sum_{a=1}^r\ovOmega_a \frac{dq^a}{q^a}, \quad 
  \text{where $\ovOmega_a$ is defined by }
  \Qp_a\cdot\ovPhi(v)=\ovPhi(\ovOmega_a(v)) \text{ for }v\in V.
\]
Then we have $\ovOmega_a(q,\hbar)=\Omega_a(q,-\hbar)$.    
Similarly, we can consider flat connections 
$\ov{\nabla}^0$, $\ov{\nabla}^1$ of the endomorphism bundle 
which correspond to $\nabla^0$, $\nabla^1$. 
There exists a unique flat section 
$\ov{S}(q,\hbar)$ of $\ov{\nabla}^1$ 
defined in the neighborhood of $q=0$ 
with the initial condition $\ov{S}(0,\hbar)=\id$. 
Then we have $\ov{S}(q,\hbar)=S(q,-\hbar)$. 
Put $\Xi\defequal S^{-1}\circ \Phi^{-1}$ and 
$\ov{\Xi}\defequal \ov{S}^{-1}\circ \ovPhi^{-1}$. 
We define a pairing $(\cdot,\cdot)$ between $\ovE$ and $E$ by 
\begin{equation*} 
(\cdot,\cdot)\colon 
\ovE\times E\rightarrow \C[\hbar,\hbar^{-1}]\!][\![q^1,\dots,q^r]\!], 
\qquad 
(x,y)\defequal
\left\langle \ov{\Xi}(x), \Xi(y)\right\rangle, 
\quad x\in \ovE, y\in E.
\end{equation*} 

It is easy to see that $\Xi$ is transformed as 
$\hXi=\Xi\circ e^{p(\log q-\log \hq)/\hbar+F/\hbar}$ 
by a change of frames $\Phi\mapsto \hPhi$ and generalized coordinates 
$(\Qq,\Qp)\mapsto (\hQq,\hQp)$.   
Therefore the pairing $(\cdot,\cdot)$ is independent of 
a choice of a frame and generalized coordinates. 
Because $\nabla^1S=0$, 
the map $\Xi=S^{-1}\Phi^{-1}$ 
satisfies the following differential equation. 
\begin{lemma}
\label{lem:differentialequationforXi}
\[\Xi(\Qp_a\cdot x)=(\hbar\parti{a}+p_a)\Xi(x), 
\]
\end{lemma}  
From this we deduce the Frobenius property: 
\[ \hbar\parti{a}(x,y)=(-\Qp_a\cdot x,y)+(x,\Qp_a\cdot y).\]
Because the pairing $\pair{\cdot}{\cdot}$ is symmetric, we have 
\[ \ov{(x,y)}=(\ov{y},\ov{x}). 
\]
When the pairing $\pair{\cdot}{\cdot}$ is non-degenerate, 
we may define $J_{\Phi,(\Qq,\Qp)}$ by the following relation. 
\[\pair{v}{J_{\Phi,(\Qq,\Qp)}}
  =(\ovPhi(\ov{L}_{\Phi,(\Qq,\Qp)}(v)),\Phi(e_0)),\quad 
  v\in V.            
\]
This corresponds to the definition of the $J$-function by Givental. 

Further, suppose that $E$ is a graded quantum $D$-module 
and that there exists a constant $n$ such that 
\[\pair{\mu v}{w}+\pair{v}{\mu w}=n\pair{v}{w}. 
\]
The original quantum $D$-module arising from 
the quantum cohomology of $M$ satisfies
the above equation for $n=\dim_{\R}M$. 
In this case, the pairing $(\cdot,\cdot)$ satisfies 
\[(\Eul+n)(x,y)=(\adE(x),y)+(x,\adE(y)). 
\]
This follows from $[\adE,\Xi]=0$. 
(Here, $\adE$ acts on $\ov{E}$ as $\adE(\ov{x})=\ov{\adE(x)}$. )

Finally, we discuss the polynomiality of the pairing in $\hbar$. 
In the case of original quantum $D$-modules, 
we have a polynomial pairing with respect to $\hbar$.  
(In fact, it coincides with the Poincar\'{e} pairing.)
In general, however, 
the pairing $(\cdot,\cdot)$ 
does not necessarily take its values in $\C[\hbar][\![q]\!]$ 
and may contain negative powers of $\hbar$.  
For a canonical frame $\Phi_{\rm can}$, we have 
$S_{\rm can}^{-1}=\id+O(\hbar^{-1})$, which we showed 
before equation (\ref{eq:propertyofcanonicalJ}). 
Hence we have
\begin{equation}
\label{eq:pairingandcanonicalframe}
(\ovPhi_{\rm can}(v),\Phi_{\rm can}(w))=\pair{v}{w}
      +O(1/\hbar), \text{ for }v,w \in V. 
\end{equation}

\begin{proposition}[polynomiality]
\label{prop:polynomiality}
The following statements are equivalent. 

$(1)$ $\pair{v}{\Omega_a^{\rm can}(w)}=
       \pair{\Omega_a^{\rm can}(v)}{w}$ 
       for $a=1,\dots,r$, $v,w\in V$.   

$(2)$ $(\ovPhi_{\rm can}(v),\Phi_{\rm can}(w))
      =\pair{v}{w}$ for $v,w\in V$.  

$(3)$ $\pair{\ov{S}(v)}{S(w)}=\pair{v}{w}$ for 
      $v, w\in V$, where $S=S_{\Phi_{\rm can},(\Qq,\Qp)}$. 

$(4)$ The pairing $(\cdot,\cdot)$ takes values in $\C[\hbar][\![q]\!]$. 

$(5)$ The function $\mathcal{G}(q,\hbar,z)\defequal 
      \pair{J_{\Phi,(\Qq,\Qp)}(q,-\hbar)}
      {J_{\Phi,(\Qq,\Qp)}(qe^{\hbar z},\hbar)}$
      is polynomial in $\hbar$, 
      where $pz\defequal\sum_{a=1}^r p_az^a$ and 
      $qe^{\hbar z}=\{q^ae^{\hbar z^a}\}_{a=1}^r$.   
\end{proposition}           
\begin{proof}
$(2)\Leftrightarrow (3)$: This follows from the definition.   

$(1)\Leftrightarrow (3)$:  
Suppose $(1)$ holds.  
Then we have 
\begin{eqnarray*}
\hbar\parti{a}\pair{\ov{S}(v)}{S(w)}&=&
   \pair{\Omega_a^{\rm can}\ov{S}(v)-\ov{S}(p_a v)}{S(w)} \\
   & &\qquad \qquad \quad
   +\pair{\ov{S}(v)}{-\Omega_a^{\rm can}{S}(w)+S(p_a w)} \\
   &=& \pair{-\ov{S}(p_a v)}{S(w)}+
    \pair{\ov{S}(v)}{S(p_a w)} 
\end{eqnarray*}
by the assumption. 
Set $F(v,w)\defequal \sum_{d\ge 0}F_d(v,w)q^d
\defequal \pair{\ov{S}(v)}{S(w)}$. 
Denote by $\ad(p_a)$ the operator acting on bilinear forms 
defined by 
\[\ad(p_a)\colon B(v,w)\longmapsto -B(p_a v,w)+B(v,p_a w), \]
where $B(\cdot,\cdot)$ is a bilinear form with values in 
$\C[\hbar,\hbar^{-1}]\!]$. 
Then, by the above calculation, 
$(\hbar d_a-\ad(p_a))F_d=0$ holds.  
Because the operator $\hbar d_a-\ad(p_a)$ is invertible 
for $d_a\neq 0$, we have 
$F(v,w)=F_0(v,w)=\pair{v}{w}$. 
Conversely, suppose $\pair{\ov{S}(v)}{S(w)}=\pair{v}{w}$ 
holds for all $v,w\in V$. 
Then, 
\begin{eqnarray*}
0&=&\hbar\parti{a}\pair{v}{w}
 =\hbar\parti{a}\pair{\ov{S}(v)}{S(w)} \\
 &=&\pair{\Omega_a^{\rm can}\ov{S}(v)}{S(w)}
  -\pair{\ov{S}(v)}{\Omega_a^{\rm can}S(w)}
  -\pair{p_a v}{w}+\pair{v}{ p_a w} \\
 &=&\pair{\Omega_a^{\rm can}\ov{S}(v)}{S(w)}
  -\pair{\ov{S}(v)}{\Omega_a^{\rm can}S(w)}.
\end{eqnarray*} 

$(2)\Leftrightarrow (4)$: This follows from the equation 
(\ref{eq:pairingandcanonicalframe}).

$(4)\Leftrightarrow (5)$: Note that $(3)$ is equivalent to the fact 
that $(\ovPhi(e_0),e^{\Qp z}\Phi(e_0))$ is polynomial in $\hbar$. 
By Lemma \ref{lem:differentialequationforXi} and 
$J_{\Phi,(\Qq,\Qp)}(q,\hbar)=e^{p\log q/\hbar}\Xi(\Phi(e_0))$, 
we have 
\begin{eqnarray*}
\Xi(e^{\Qp z}\Phi(e_0))
 =e^{\sum_{a=1}^r z^a\left(\hbar\lparfrac{}{q^a}+p_a\right)}\Xi(\Phi(e_0))
 =e^{-{p\log q}/\hbar}J_{\Phi,(\Qq,\Qp)}(qe^{\hbar z},\hbar).
\end{eqnarray*} 
Therefore we have $\mathcal{G}(q,\hbar,z)=
\langle\ovXi(\Phi(e_0)),\Xi(e^{\Qp z}\Phi(e_0))\rangle =
(\ovPhi(e_0),e^{\Qp z}\Phi(e_0))$. 
\end{proof}

\begin{remark}
\label{rem:pairing} 
When an abstract quantum $D$-module $E$ 
arises from quantum cohomology, 
we take $\pair{\cdot}{\cdot}$ to be the Poincar\'{e} pairing  
on $H^*(M,\C)$. 
Then by equation (\ref{eq:Frobenius}) 
in Section \ref{sect:quantumcohomologyandDmodule},  
the pairing $(\cdot,\cdot)$ satisfies the polynomiality 
and we have $(x,y)=\pair{x}{y}$ ($\Phi=\id$ in this case).  
(5) is the famous polynomiality property of the $J$-function. 
It imposes a strong constraint on the $J$-function and 
Givental used it to calculate the $J$-function 
\cite{givental-mirrorthm-projective,givental-mirrorthm-toric}. 
\end{remark}


\section{Equivariant Floer Theory} 
\label{sect:equivariantFloertheory}

In this section, we explicitly construct 
equivariant Floer cohomology as an abstract quantum $D$-module 
for toric complete intersections. 
More precisely, we construct it for a superspace 
$(\XSig,\mV)$ where $\XSig$ is a toric manifold and 
$\mV$ is a sum of nef line bundles. 
In the case where $c_1(\XSig/\mV)\defequal c_1(T\XSig)-c_1(\mV)$ 
is nef, we prove that 
the equivariant Floer cohomology is isomorphic to 
the quantum $D$-module as an abstract quantum $D$-module.   

\subsection{Toric Manifolds}

To fix the notation, we briefly recall 
the definition of toric manifolds following 
Audin's book \cite{audin}. 
Let $\Sigma$ be a fan in $\R^n$. 
Let $\Sigma^{(1)}=\{x_1,\dots,x_N\}$ 
be its $1$-skeleton, where 
the $x_i$'s are primitive vectors in $\Z^n$ 
which denote generators of one dimensional cones. 
We set $r\defequal N-n$. 
For a subset $I$ of $\{1,2,\dots,N\}$, we define  
$\sigma_{I}$ to be the cone generated by $\{x_i\}_{i\in I}$. 
We assume that $\Sigma$ is smooth and complete 
in the sense that  
\begin{eqnarray*}
\Sigma\text{ is smooth} &\Longleftrightarrow& 
\text{each cone }\sigma \in\Sigma
\text{ is generated by part of a $\Z$ basis} \\
\Sigma\text{ is complete} &\Longleftrightarrow& 
\text{its support } |\Sigma|=\cup_{\sigma\in \Sigma}\sigma 
\text{ is equal to }\R^n. 
\end{eqnarray*} 
Let $\pi$ be the following homomorphism. 
\[\pi\colon\Z^N\rightarrow \Z^n,\qquad
e_i\mapsto x_i,\]
where $\{e_i\}_{i=1}^N$ is a canonical basis of $\Z^N$. 
We extend the homomorphism $\pi$ to 
$\pi_{\C}=\pi\otimes\id_{\C}\colon \C^N\rightarrow \C^n$. 
It satisfies $\pi_{\C}(\Z^N)\subset \Z^n$, therefore 
it induces a homomorphism $\bar{\pi}_{\C}$ of tori 
$\bar{\pi}_{\C}\colon \C^N/\Z^N \rightarrow \C^n/\Z^n$. 
We define  
\[K_{\C}\defequal
\ker(\bar{\pi}_{\C}\colon T_{\C}^N\rightarrow T_{\C}^n).\]
In the same way, we define $\pi_{\R}\colon \R^N\rightarrow \R^n$ 
and $K_{\R}\defequal \ker(\bar{\pi}_{\R}\colon\R^N/\Z^N\rightarrow 
\R^n/\Z^n)$. 
Because $\Sigma$ is complete and smooth, 
$K_{\C}$ turns out to be connected and isomorphic to $T_{\C}^r$. 
$K_{\C}$ is a subgroup of $\C^N/\Z^N\cong T_{\C}^N$. 
Let $(t_1,\dots,t_N)$ be the coordinates of $T_{\C}^N$. 
We have an exact sequence. 
\begin{equation}
\label{eq:exactseq}  
\begin{CD}
0 @>>> \Lie(K_{\C}) @>{i}>> \C^N @>{\pi_{\C}}>> \C^n @>>> 0.
\end{CD}
\end{equation} 
For a subset $I$ of $\{1,2,\dots,N\}$, 
let $\ov{I}$ denote the complement of $I$ and 
$\C^I$ denote the vector space 
\[\C^I=\{(z_1,\dots,z_N)\in \C^N\;
|\;i\notin I\Longrightarrow z_i=0\}. \]
Similarly, for any set $Y$, 
we define $Y^I$ as a subset of $Y^N$.  
Define an open set $U_{\Sigma}$ in $\C^N$ as    
\[U_{\Sigma}\defequal \C^N-
\bigcup_{\sigma_{I}\notin\Sigma}\C^{\ov{I}}. \] 
The torus $K_{\C}$ acts on $U_{\Sigma}$ freely as a subgroup of 
$T_{\C}^N$. 
\[U_{\Sigma}\ni(z_1,\dots,z_N)\mapsto (t_1z_1,\dots,t_Nz_N) 
\quad \text{for } (t_1,\dots,t_N)\in K_{\C}. \] 
We define the {\bf toric manifold} $\XSig$ as the quotient  
\[\XSig\defequal U_{\Sigma}/K_{\C}.\]
It turns out that $\XSig$ is smooth and compact, and that  
$U_{\Sigma}$ is a principal $K_{\C}$ bundle 
over $\XSig$. 

The Picard group $\Pic(\XSig)$ 
is identified with $\Hom(K_{\C},T_{\C})$ 
and also with $H^2(\XSig,\Z)$ as follows:  
\begin{eqnarray*}
 \Z^r\cong\Hom(K_{\C},T_{\C})\ni \rho 
&\mapsto& \mathcal{L}_{\rho}\defequal 
\C\times U_{\Sigma}/(\rho(t)^{-1}v,z)\sim(v,tz),\,t\in K_{\C}
\;\in \Pic(\XSig) \\
&\mapsto& c_1(\mathcal{L}_\rho)\in H^2(\XSig,\Z). 
\end{eqnarray*} 
Let $e^*_i$ be the character of $T_\C^N$ defined by     
$e_i^*(t_1,\dots,t_N)\defequal t_i$. 
Composing it with the inclusion $K_{\C}\hookrightarrow T_{\C}^N$, 
we have a character of $K_{\C}$ 
and the corresponding Chern class $u_i$ in $H^2(\XSig,\Z)$. 
The anti-canonical class can be written as 
$c_1(T\XSig)=u_1+u_2+\dots+u_N$. 
The abelian group $\Hom(K_{\C},T_{\C})$ is an integral 
lattice of $\Lie(K_{\C})^\vee$
and we have the following identifications. 
\begin{alignat*}{3} 
H^2(\XSig,\R)&\cong \Lie(K_{\R})^\vee,&\quad 
H^2(\XSig,\C)&\cong \Lie(K_{\C})^\vee,&\quad 
H^2(\XSig,\Z)&\cong \Hom(K_{\C},T_{\C}),\\
H_2(\XSig,\R)&\cong \Lie(K_{\R}),&\quad 
H_2(\XSig,\C)&\cong \Lie(K_{\C}),&\quad
H_2(\XSig,\Z)&\cong \Hom(T_{\C},K_{\C}). 
\end{alignat*}

We can describe the K\"{a}hler cone of $H^2(\XSig,\R)
=H^{1,1}(\XSig,\R)$.  
Let $M\colon (\R^N)^\vee\rightarrow \Lie(K_{\R})^\vee$ 
denote the dual of the inclusion 
$i\colon \Lie(K_{\R})\hookrightarrow \R^N$.   
Then we have 
\begin{equation}
\label{eq:Kahlercone}
\text{K\"{a}hler cone }\defequal \{\text{nef classes in 
$H^2(\XSig,\R)$}\}
=\bigcap_{\sigma_{I}\in \Sigma}
M(\R_{+}^{\ov{I}}), 
\end{equation}
where $\R_{+}$ denotes the set of non-negative real numbers. 
We assume that 
the interior of $\bigcap_{\sigma_{I}\in \Sigma}M(\R_{+}^{\ov{I}})$ 
is non-empty so that $\XSig$ becomes a K\"{a}hler manifold. 
We can describe $\XSig$ as a symplectic reduction of $U_\Sigma$.  
Let $\mu\colon U_{\Sigma}\rightarrow \Lie(K_{\R})^{\vee}$ be a moment map 
of the $T^r_{\R}$ action. 
\[\mu(z_1,\dots,z_N)=M(\frac{1}{2}|z_1|^2,\dots,\frac{1}{2}|z_N|^2).\] 
We choose an element $\eta$ 
in the interior of the K\"{a}hler cone. 
Then we have an isomorphism 
\[\XSig\cong \mu^{-1}(\eta)/T_{\R}^r.\]
The reduced symplectic form gives a K\"{a}hler form of 
$\XSig$ and represents the cohomology class $\eta$. 

Let $\mV=\mV_1\oplus \dots \oplus\mV_l$ 
be a sum of line bundles over $\XSig$ 
such that each $v_i\defequal c_1(\mV_i)$ is nef. 
Then we have a convex superspace $(\XSig,\mV)$. 
As in Section \ref{sect:quantumcohomologyandDmodule}, 
we take nef classes $p_1,\dots,p_r$ such that they  
form a basis of $H^2(\XSig,\Z)$. 
Moreover, when $c_1(\XSig/\mV)$ is nef, 
we choose the $p_i$'s so that $c_1(\XSig/\mV)$ is contained in 
the cone generated by $p_1,\dots,p_r$. 
We put $u_i=\sum_{a=1}^rm_i^a p_a$ and $v_j=\sum_{a=1}^r l_j^a p_a$. 
The matrix $(m_i^a)$ represents $M$ with respect to
the basis $\{e^*_1,\dots,e^*_N\}$ of $(\C^N)^\vee$ and 
the basis $\{p_1,\dots,p_r\}$ of $\Lie(K_{\C})^\vee$   
because $u_i=M(e_i^*)$. 
We establish the following lemma for later use. 
\begin{lemma}
\label{lem:Moricone}
Let $\Lambda$ be the semigroup in $H_2(\XSig,\Z)$ 
generated by effective curves in $\XSig$. 
If $d\notin \Lambda$, we have 
$\sigma_{I_0}\notin \Sigma$ for 
$I_0\defequal\{i\in \{1,\dots,N\}|\;\pair{u_i}{d}<0\}$. 
\end{lemma} 
\begin{proof} 
If $\sigma_{I_0}\in\Sigma$, we have $\pair{x}{d}\ge 0$ 
for any $x\in M(\R^{\ov{I}_0}_{+})
=\{\sum_{i\notin I_0}c_iu_i|\;c_i\ge 0\}$. 
By equation (\ref{eq:Kahlercone}), 
we conclude that $x$ is in the dual cone of the K\"{a}hler cone which 
is equal to $\Lambda$.  
\end{proof}

The cohomology ring of the toric manifold $\XSig$ is 
generated by $p_1,\dots,p_r$ and subject to the following relations, 
see e.g. \cite{brion}:
\begin{equation}
\label{eq:relationoftoric}
u_{i_1}u_{i_2}\dots u_{i_l}=0, \quad 
\text{if }\sigma_{\{i_1,\dots,i_l\}}\notin \Sigma. 
\end{equation}

\subsection{Givental's Model for Free Loop Spaces}

We describe Givental's model 
for the universal covering of the free loop space of $\XSig$. 
Originally, Givental considered this model for $\Proj^n$ in 
\cite{givental-homologicalgeometry}. 
The model for toric manifolds was given by 
Vlassopoulos \cite{vlassopoulos}. 

In order to understand the model, 
we consider a more general situation. 
Let $X$ be a simply-connected manifold and 
$\{p_1,\dots,p_r\}$ be a basis of $H^2(X,\Z)$. 
Let $\mathcal{L}(p_i)$ be a line bundle whose first Chern class 
is $p_i$. 
We set $U=\prod_{i=1}^r(\mathcal{L}(p_i)\setminus \{\text{zero section}\})$ 
which is a principal $T_{\C}^r$ bundle over $X$. 
It is easy to see that $U$ is $2$-connected. 
In the case of toric manifolds, $X$ corresponds to 
$\XSig$ and $U$ corresponds to $U_{\Sigma}$. 
There exists a homeomorphism 
\begin{equation}
\label{eq:universalcover}
LU/L_0T_{\C}^r\cong \widetilde{LX},
\end{equation} 
where $LU$ denotes the free loop space of $U$, 
$L_0T_{\C}^r$ denotes the set of 
all contractible loops in $T_{\C}^r$ and  
$L_0T_{\C}^r$ acts on $LU$ 
by pointwise multiplication. 
Given a loop $\tilde{\gamma}\colon S^1\rightarrow U$, 
we have a disk $\tilde{g}\colon D^2\rightarrow U$ 
which contracts the loop $\tilde{\gamma}$ 
i.e. $\tilde{g}|_{\partial D^2}=\tilde{\gamma}$. 
This disk $\tilde{g}$ is uniquely determined up to homotopy 
because $U$ is $2$-connected. 
Then, by composing with the projection $U\rightarrow X$, 
we obtain a pair $(\gamma,[g])$ consisting of 
a loop $\gamma$ in $X$ and a homotopy class of 
the contracting disk $g$ in $X$ 
which represents an element of $\widetilde{LX}$. 
Conversely, if we have an element $(\gamma,[g])$ 
of $\widetilde{LX}$, 
we can lift $g$ to the map $\tilde{g}\colon D^2\rightarrow U$ 
and then have a loop $\tilde{\gamma}=\tilde{g}|_{S^1}$ in $U$.  

Givental's model corresponds to the left hand side of 
the homeomorphism (\ref{eq:universalcover}), 
but contains  only polynomial loops. 
Let $\C[\zeta,\zeta^{-1}]$ be a Laurent polynomial ring. 
We define Givental's model $L$ as  
\begin{eqnarray*}
LU_{\Sigma}\defequal \C[\zeta,\zeta^{-1}]^N-
\bigcup_{\sigma_{I}\notin\Sigma}\C[\zeta,\zeta^{-1}]^{\ov{I}}, \qquad 
L\defequal LU_{\Sigma}/T_{\C}^r.  
\end{eqnarray*} 
Givental's model $L$ is an infinite dimensional toric variety  
which may be defined by an infinite dimensional fan. 
The model $L$ can also be written 
as a symplectic reduction of $LU_{\Sigma}$. 

Any element of $L$ can be represented by an $N-$tuple of 
Laurent polynomials $(\gamma_1(\zeta),\dots,\gamma_N(\zeta))$, 
where $\gamma_i\in\C[\zeta,\zeta^{-1}]$. 
For a generic value of $\zeta\in \C$,   
$(\gamma_1(\zeta),\dots,\gamma_N(\zeta))$ is 
contained in $U_{\Sigma}$, 
but there may be finitely many values $\zeta$ such 
that $(\gamma_1(\zeta),\dots,\gamma_N(\zeta))$ 
is not contained in $U_{\Sigma}$.   
Therefore, we can regard $[\gamma_1,\dots,\gamma_N]$ as a 
holomorphic map 
\[[\gamma_1,\dots,\gamma_N]\colon 
\Proj^1\setminus \{\text{finite points}\}
\rightarrow \XSig,\]
where $\Proj^1=\C\cup\{\infty\}$. 
Because $\XSig$ is a complete variety, 
we can extend this map on the whole $\Proj^1$. 
Consequently we have a (not continuous) map 
\[L\rightarrow \Hol(\Proj^1,\XSig),\]
where $\Hol(\Proj^1,\XSig)$ denotes the set of holomorphic maps 
from $\Proj^1$ to $\XSig$. 
However, some information is lost when we pass from $L$
to $\Hol(\Proj^1,\XSig)$. 
For example, take a tuple 
$(t_1(\zeta),\dots,t_N(\zeta))$ in  
$\C[\zeta,\zeta^{-1}]^N$ such that 
$t_i(\zeta)\neq 0$ and 
for a generic value $\zeta\in\C$, 
$(t_1(\zeta),\dots,t_N(\zeta))$ is contained in $K_{\C}$. 
Then two elements 
$[\gamma_1,\dots,\gamma_N]$ and 
$[t_1\gamma_1,\dots,t_N\gamma_N]$ 
give the same element in $\Hol(\Proj^1,\XSig)$.  

We define an $S^1$ action on the model $L$ by  
\[(\gamma_1(\zeta),\dots,\gamma_N(\zeta))\mapsto 
(\gamma_1(e^{\sqrt{-1}\theta}\zeta),\dots,
          \gamma_{N}(e^{\sqrt{-1}\theta}\zeta)),
\quad e^{\sqrt{-1}\theta}\in S^1.\] 
This action is Hamiltonian with respect to the   
K\"{a}hler forms obtained by 
symplectic reduction. 
The Hamiltonian $H$ can be written as  
\[H([\gamma_1,\dots,\gamma_N])=\sum_{i=1}^N\sum_{\nu=-\infty}^\infty
                             \frac{1}{2}\nu|a_{i,\nu}|^2, 
 \quad  \gamma_i(\zeta)=\sum_{\nu}a_{i,\nu}\zeta^\nu\]
on the level set of the symplectic reduction. 
The $S^1$-fixed components are all isomorphic to $\XSig$ 
and each of them is parameterized by $H_2(\XSig,\Z)$. 
For $d\in H_2(\XSig,\Z)$, we define 
\[(\XSig)_d\defequal 
\left\{(z_1\zeta^{\pair{u_1}{d}},\dots,z_N\zeta^{\pair{u_N}{d}})
       \in L\right\}. \]
This is fixed by $S^1$ because  
$\{(e^{\sqrt{-1}\pair{u_1}{d}\theta},\dots,
   e^{\sqrt{-1}\pair{u_N}{d}\theta})\}$ 
forms a one-parameter subgroup of $K_{\C}$. 
In this paper, Floer theory means Bott-Morse theory 
on the infinite dimensional variety $L$.  
We take the Hamiltonian $H$ as the Bott-Morse function. 
The gradient flow of $H$ can be written as 
\[\phi_t([\gamma_1(\zeta),\dots,\gamma_N(\zeta)])
  =[\gamma_1(\zeta e^{-t}),\dots,\gamma_N(\zeta e^{-t})], 
\quad t\in \R.
\]
The critical set of the Morse function $H$ is identical with the 
$S^1$ fixed point set $\coprod_d (\XSig)_d$. 

Givental's model has a stratification by  
finite dimensional toric manifolds. 
Take $d_1$, $d_2$ in  $H_2(\XSig,\Z)$ and define 
\[L_{d_1}^{d_2}\defequal
 \left\{(\gamma_1(\zeta),\dots,\gamma_N(\zeta))\in L 
       \Bigg| \gamma_i(\zeta)=
       \begin{cases}
            \sum_{\nu=\pair{u_i}{d_1}}^{\pair{u_i}{d_2}}
            a_{i\nu}\zeta^\nu \; \text{if } 
            \pair{u_i}{d_2-d_1}\ge 0, \\
            0 \; \; \qquad \;\;\text{otherwise}
       \end{cases}
 \right\} 
\] 
Note that $L_{d_1}^{d_2}$ is empty 
if $\sigma_{I_0}\notin \Sigma$ for the subset 
$I_0\defequal \{i\;|\;\pair{u_i}{d_2-d_1}<0\}$, and  
$L_{d_1}^{d_2}$ is the closure of the union of all gradient trajectories 
which connect two critical submanifolds 
$(\XSig)_{d_1}$ and $(\XSig)_{d_2}$: 
\[L_{d_1}^{d_2}=\overline{
   \{\gamma\in L|\lim_{t\rightarrow \infty}\phi_t(\gamma)\in (\XSig)_{d_1}, 
      \;        \lim_{t\rightarrow -\infty}\phi_t(\gamma)\in (\XSig)_{d_2}
    \}                    }. 
\]
We first claim that the union of all these spaces is equal to $L$. 
\[\bigcup_{d_1,d_2}L_{d_1}^{d_2}=L. \]
To see this, it suffices to show that 
there exists  $d\in H_2(\XSig,\Z)$ such that $\pair{u_i}{d}>0$ for 
all $i$. 
Because $\Sigma$ is complete, we have positive integers 
$r_1,\dots, r_N$ such that 
$\pi(\sum_{i=1}^Nr_ie_i)=\sum_{i=1}^Nr_ix_i=0$. 
By the exact sequence (\ref{eq:exactseq}), 
there exists a $d\in H_2(\XSig,\Z)$ such that 
$i(d)=\sum_{i=1}^Nr_ie_i$. 
Then we have $\pair{u_i}{d}=\pair{e^*_i}{\sum_{i=1}^Nr_ie_i}=r_i>0$. 

These finite dimensional spaces were introduced by Givental 
\cite{givental-mirrorthm-toric} 
under the name ``toric map space''. 
They are also toric manifolds which are smooth and compact. 
$L_{d_1}^{d_2}$ can be considered as a 
compactification of 
the space of degree $d_2-d_1$ maps from $\Proj^1$ to $\XSig$. 
In fact, when $I_0=\emptyset$, generic elements of it 
represent degree $d_2-d_1$ maps from 
$\Proj^1$ to $\XSig$. 
The complex dimension of $L_{d_1}^{d_2}$ is equal to 
$\sum_{i\notin I_0} \pair{u_i}{d_2-d_1}+n-|I_0|$. 
On the other hand, the virtual dimension of the 
space of degree $d_2-d_1$ maps is equal to 
$\pair{c_1(M)}{d_2-d_1}+\dim\XSig
=\sum_{i=1}^N \pair{u_i}{d_2-d_1}+n$. 
Therefore, we have the inequality 
\[\dim L_{d_1}^{d_2} \ge 
  \virtualdim \Hol_{d_2-d_1}(\Proj^1,\XSig). 
\]
In the case $I_0=\emptyset$, we have equality. 
The reason why the dimension is sometimes larger than 
the virtual one is that 
the stable manifold $L_{d_1}^\infty$ and 
the unstable manifold $L_{-\infty}^{d_2}$ 
do not necessarily intersect transversely, 
where 
\begin{eqnarray*}
L_{d_1}^{\infty}\defequal
 \left\{(\gamma_1(\zeta),\dots,\gamma_N(\zeta))\in L 
       \Bigg| \gamma_i(\zeta)=
              \sum_{\nu\ge\pair{u_i}{d_1}}
              a_{i\nu}\zeta^\nu
 \right\}, \\
L_{-\infty}^{d_2}\defequal
 \left\{(\gamma_1(\zeta),\dots,\gamma_N(\zeta))\in L 
       \Bigg| \gamma_i(\zeta)=
              \sum_{\nu\le \pair{u_i}{d_2}}
              a_{i\nu}\zeta^\nu
 \right\}.  
\end{eqnarray*} 
We need to deal with these intersections not geometrically 
but cohomologically 
so that we have the `correct' intersection theory.   
For this, we later define 
the Poincar\'{e} duals of semi-infinite dimensional 
spaces $L_{-\infty}^{d_2}$ and 
$L_{d_1}^{\infty}$. 

Define $S^1$ equivariant line bundles 
$\mathcal{L}_{i,\nu}$ and $\mathcal{V}_{i,\nu}$ over $L$ by 
\begin{align*}
\mathcal{L}_{i,\nu}\defequal &\C\times LU_{\Sigma}/
(t_i^{-1}v,(a_{j\mu})_{j,\mu})\sim(v,(t_ja_{j\mu})_{j,\mu}),\;
t\in K_{\C}, \\
&[v,(a_{j\mu})_{j,\mu}]\mapsto 
[e^{\sqrt{-1}\nu\theta} v,(e^{\sqrt{-1}\mu\theta} a_{j\mu})_{j,\mu}], \quad 
 \text{ for } e^{\sqrt{-1}\theta}\in S^1, \\
\mathcal{V}_{i,\nu}\defequal &\C\times LU_{\Sigma}/
(v_i(t)^{-1}v,(a_{j\mu})_{j,\mu})\sim(v,(t_ja_{j,\mu})_{j,\mu}), \;
t\in K_{\C}, \\
&[v,(a_{j\mu})_{j,\mu}]\mapsto 
[e^{\sqrt{-1}\nu\theta} v,(e^{\sqrt{-1}\mu\theta} a_{j\mu})_{j,\mu}], \quad 
 \text{ for } e^{\sqrt{-1}\theta}\in S^1, 
\end{align*} 
where $v_i(t)$ denotes the character corresponding to the class 
$v_i=c_1(\mV_i)$.  
The line bundle $\mathcal{L}_{i,\nu}$ 
has the following equivariant section $s_{i,\nu}$. 
\[s_{i,\nu}\colon [(a_{j,\mu})_{j,\mu}] \mapsto 
     [a_{i,\nu},(a_{j,\mu})_{j,\mu}].\]
The stable manifold $L_{d_1}^{\infty}$ can be described 
in terms of the sections $s_{i,\nu}$ as 
\[L_{d_1}^{\infty}=\text{zero-locus of the section }
  \prod_{j=1}^N\left(\prod_{\mu< \pair{u_j}{d_1}} s_{j,\mu}\right).
\]
We also define $S^1$ equivariant line bundle $\mathcal{L}(p_a)$ 
over $L$ by  
\begin{gather*}
\mathcal{L}(p_a)\defequal \C\times LU_{\Sigma}/
(p_a(t)^{-1}v,(a_{j\mu})_{j,\mu})\sim(v,(t_j a_{j\mu})_{j,\mu}),\;
t\in K_{\C}, \\
 [v,(a_{j\mu})_{j,\mu}]\mapsto  
 [v,(e^{\sqrt{-1}\mu\theta} a_{j\mu})_{j,\mu}], \quad 
\text{ for } e^{\sqrt{-1}\theta}\in S^1,  
\end{gather*} 
where $p_a(t)$ denotes the character corresponding to 
the class $p_a$. 
Introduce another $S^1$ action on $\mathcal{V}_{i,\nu}$ 
as scalar multiplication on its fibers. 
This $S^1$ acts on $\mathcal{L}_{i,\nu}$ and $\mathcal{L}(p_a)$ trivially. 
Thus we have $T^2$ equivariant line bundles $\mathcal{L}_{i,\nu}$, 
$\mathcal{V}_{i,\nu}$ and $\mathcal{L}(p_a)$. 
Let $P_a\defequal c_1^{T^2}(\mathcal{L}(p_a))$, 
$U_{i,\nu}\defequal c_1^{T^2}(\mathcal{L}_{i,\nu})$ and 
$V_{i,\nu}\defequal c_1^{T^2}(\mathcal{V}_{i,\nu})$ 
denote the equivariant first Chern classes. 
Let $\hbar,\lambda$ denote generators of 
the $T^2$ equivariant cohomology of a point, 
where $\hbar$ is a generator of the $S^1$-action rotating loops 
and $\lambda$ is that of the fiberwise $S^1$-action. 
More precisely, let $\C w^n {w'}^m$ denote a rank one $T^2$ module  
on which $T^2$ acts as 
$w^n {w'}^m\mapsto e^{\sqrt{-1}(n\theta+m\theta')} w^n{w'}^m$,
$(e^{\sqrt{-1}\theta},e^{\sqrt{-1}\theta'})\in T^2$. 
Then we define 
$\hbar\defequal c_1^{T^2}(\C w^{-1}\rightarrow \mathrm{pt})$ 
and $\lambda\defequal c_1^{T^2}(\C {w'}^{-1}\rightarrow \mathrm{pt})$. 

\begin{lemma} 
\label{lem:UVPrelation}
The equivariant Chern classes $U_{i,\nu}$, $V_{i,\nu}$ can be  
written as 
\[ U_{i,\nu}=\sum_{a=1}^rm_i^aP_a-\nu\hbar,\qquad 
   V_{i,\nu}=\sum_{a=1}^r l_i^aP_a-\nu\hbar-\lambda. 
\]
The restriction of $P_a$ to the fixed component gives  
\[P_a|_{(\XSig)_d}=p_a+\pair{p_a}{d}\hbar.
\] 
\end{lemma}
\begin{proof} 
The proposition follows from the isomorphisms 
$\mathcal{L}_{i,\nu}\cong 
\bigotimes_{a=1}^r\mathcal{L}(p_a)^{\otimes m_i^a}\otimes \C w^{\nu}$, 
$\mathcal{V}_{i,\nu}\cong
\bigotimes_{a=1}^r\mathcal{L}(p_a)^{\otimes l_i^a}\otimes 
\C w^{\nu}{w'}^{-1}$ and
$L(p_a)|_{(\XSig)_d}\cong 
L(p_a)|_{(\XSig)_0}\otimes \C w^{-\pair{p_a}{d}}$ 
as $T^2$ equivariant line bundles. 
\end{proof}

Next we define covering transformations on $L$. 
As in Section \ref{sect:quantumcohomologyandDmodule}, 
let $\{q^1,\dots,q^r\}$ denote a basis of $H_2(\XSig,\Z)$ dual to 
$\{p_1,\dots,p_r\}$. 
Each $q^a$ defines the following covering transformation 
$Q^a\colon L\rightarrow L$.
\[Q^a\colon (\gamma_1(\zeta),\dots,\gamma_N(\zeta))
  \mapsto(\zeta^{-m_1^a}\gamma_1(\zeta),\dots,\zeta^{-m_N^a}\gamma_N(\zeta)).
\]
Note that $m_i^a=\pair{M(e_i^*)}{q^a}=\pair{e^*_i}{i(q_a)}$ and 
$(\zeta^{m_1^a},\dots,\zeta^{m_N^a})$ is a one-parameter subgroup 
of $K_{\C}$ corresponding to $q^a$. 
Note that $Q^a$ is an $S^1$ equivariant map.
Over the finite dimensional strata, $Q^a$ defines a map 
$Q^a\colon L_{d_1}^{d_2}\rightarrow 
           L_{d_1-q^a}^{d_2-q^a}$. 
\begin{proposition}
\label{prop:commutationrelation}
The pull back of $P_a$ by $Q^b$ is 
\[{Q^b}^*(P_a)=P_a-\hbar\delta_a^b\]
Therefore we have the commutation relation 
$[P_a,(Q^b)^*]=\delta_a^b (Q^b)^*$. 
\end{proposition}
\begin{proof}
This follows from the isomorphism 
${Q^a}^*(L(p_b))\cong L(p_b)\otimes\C w^{\delta_b^a}$.  
\end{proof}

\subsection{Equivariant Floer theory} 

As already remarked, equivariant Floer theory considered in this paper 
is Bott-Morse theory on $L$ 
which uses the Hamiltonian $H$ as Bott-Morse function.  
The gradient flows and critical submanifolds were explicitly described in 
the previous subsection. 
None of the moduli spaces $L_{d_1}^{d_2}$ of gradient trajectories 
is real one-dimensional. 
Hence all the differentials in the Morse-Witten-Floer complex are zero 
and we have an isomorphism of modules:   
\[FH^*(L;H)\cong H^*(\XSig)\otimes\C[\![H_2(\XSig,\Z)]\!]. \]
We furthermore investigate a $D$-module structure on the Floer cohomology.  
For that, we realize classes in Floer cohomology 
as classes of semi-infinite degree. 
 
We formulate equivariant Floer theory in the category of superspaces. 
An object of the category is a pair consisting of a 
topological space and a (possibly infinite dimensional) vector bundle on it. 
A morphism between two objects $(M,\mV)$ and $(N,\mW)$ 
is a pair consisting of a continuous map $f\colon M\rightarrow N$ 
and a bundle map $\phi\colon f^*(\mW)\rightarrow \mV$. 
The cohomology of the superspace $(M,\mV)$ is defined by 
\[H^*(M/\mV)\defequal H^*(M).
\]
The functor $H^*$ is a contravariant functor 
from the category of superspaces to the category of graded rings. 
For a morphism $(f,\phi)\colon (M,\mV)\rightarrow (N,\mW)$, 
we have a pull back
\[(f,\phi)^*\colon H^*(N/\mW)\rightarrow H^*(M/\mV)
\]  
which coincides with the ordinary pull back by $f$. 
When $f$ is a proper map between manifolds $M$ and $N$,
$\phi$ is injective and $\cok(\phi)$ is of finite rank,   
we can define a push-forward 
\[(f,\phi)_*\colon H^*(M/\mV)\rightarrow H^*(N/\mW), \quad
\alpha\mapsto f_*(\alpha\cup Euler(\cok(\phi))).
\] 
This raises the degree by $(\dim N-\rank\mW)-(\dim M-\rank\mV)$. 
This definition of push-forward comes from the following 
definition of the fundamental class of the superspaces: 
$[M/\mV]\defequal [M]\cap Euler(\mV)$ when $M$ is a compact manifold.  
A cohomology class of the superspace $(M,\mV)$ behaves like 
the class restricted to the zero locus of a transverse section of $\mV$. 

In order to define Floer cohomology, we use the 
$T^2$-equivariant cohomology of superspaces. 
We have $T^2$-equivariant infinite dimensional vector bundles 
$\mW_{d_1}^{d_2}$ over $L$ indexed by $d_1$ and $d_2$ in $H_2(\XSig,\Z)$. 
\[\mW_{d_1}^{d_2}\defequal 
   \bigoplus_{i=1}^N\bigoplus_{\nu=-\infty}^{\pair{u_i}{d_1}-1}
   \mathcal{L}_{i,\nu} \oplus 
   \bigoplus_{i=1}^l\bigoplus_{\nu=\pair{v_i}{d_2}}^{\infty} 
   \mathcal{V}_{i,\nu}. 
\]
Because the cohomology of Givental's model $L$ 
is the polynomial ring 
generated by $P_a$, $\hbar$ and $\lambda$, we have
\footnote{In general, the $S^1$ equivariant cohomology 
of the universal cover $\widetilde{LM}$ 
of free loop space for simply-connected manifold $M$ contains 
a polynomial ring generated by two-dimensional classes of $M$. 
Furthermore, after the localization, we have an isomorphism 
\[H^*_{S^1}(\widetilde{LM},\R)\otimes_{\R[\hbar]}\R[\hbar,\hbar^{-1}]
       \cong \R[p_1,\dots,p_r,\hbar,\hbar^{-1}],\]
where $\{p_a\}_{a=1}^r$ is a linear basis of $H^2(M)$
\cite{iritanicoh}. 
}
\[H^*_{T^2}(L/\mW_{d_1}^{d_2})=\C[P_1,\dots,P_r,\hbar,\lambda].
\]

When $\pair{u_i}{d'_1}\le \pair{u_i}{d_1}$ and 
$\pair{v_i}{d_2}\le \pair{v_i}{d'_2}$ hold for all $i$, 
we have a natural morphism 
$(\id,\iota)\colon (L,\mW_{d_1}^{d_2})\rightarrow (L,\mW_{d'_1}^{d'_2})$ 
and a push-forward. 
\[(\id,\iota)_*\colon H^*_{T^2}(L/\mW_{d_1}^{d_2})
         \longrightarrow H^*_{T^2}(L/\mW_{d'_1}^{d'_2}), 
\quad \alpha\mapsto \alpha\cup
      \prod_{i=1}^N\prod_{\nu=\pair{u_i}{d'_1}}^{\pair{u_i}{d_1}-1}U_{i,\nu}
 \cup\prod_{i=1}^l\prod_{\nu=\pair{v_i}{d_2}+1}^{\pair{v_i}{d'_2}}V_{i,\nu}.
\]
This raises the degree by 
$2\sum_{i=1}^N\pair{u_i}{d_1-d'_1}+2\sum_{i=1}^l\pair{v_i}{d'_2-d_2}$. 
Therefore $\{H^*_{T^2}(L/\mW_{d_1}^{d_2})\}_{d_1,d_2}$ 
forms an inductive system and we take its limit. 
\[A^*_{T^2}\defequal \underset{d_1,d_2}{\varinjlim} H^*_{T^2}(L/\mW_{d_1}^{d_2})
  \left[\textstyle 2\sum_{i=1}^N\pair{u_i}{d_1}-2\sum_{i=1}^l\pair{v_i}{d_2}
  \right],
\]
where the bracket $[\cdots]$ means a shift of grading, 
e.g. $M[i]^n=M^{i+n}$ for a graded module $M$.  
Then $A^*_{T^2}$ becomes a graded module because of the shift of grading. 
It contains $H^*_{T^2}(L/\mW_{d_1}^{d_2})$ as a submodule. 
Each element in $A^*_{T^2}$ has an infinite product expression 
which can be considered to be a push-forward to $H^*_{T^2}(L)$ 
as follows:
\[\alpha\cup \prod_{i=1}^N\prod_{\nu<\pair{u_i}{d_1}}U_{i,\nu}
        \cup \prod_{i=1}^l\prod_{\nu\ge\pair{v_i}{d_2}}V_{i,\nu} 
  \quad \text{ for } \alpha \in H^*_{T^2}(L/\mW_{d_1}^{d_2}).
\]
This expression is compatible with the inductive limit. 

We explain a geometric meaning of elements in $A^*_{T^2}$. 
The class $1$ in $H^*_{T^2}(L/\mW_{d_1}^{d_2})$ is, 
after being pushed forward to $L$, considered to be 
the Poincar\'{e} dual of a zero-locus 
of any transverse section of $\mW_{d_1}^{d_2}$. 
Therefore it represents the fundamental class of the superspace 
\begin{equation}
\label{eq:superspaceoverstablemanifold}
(L_{d_1}^{\infty},
 \bigoplus_{i=1}^l\bigoplus_{\nu=\pair{v_i}{d_2}}^\infty 
 \mathcal{V}_{i,\nu}) 
\end{equation}
over the stable manifold $L_{d_1}^{\infty}$. 
We relate this superspace with (the model of) the loop space of the 
complete intersection in $\XSig$ when $d_1=d_2$.  
When we have a transverse section $s\colon \XSig\rightarrow \mV$ and its 
zero-locus $Y\defequal s^{-1}(0)$, 
the model $L_Y$ for the loop space of $Y$ may be written as 
\[L_Y\defequal \left\{[\gamma_1,\dots,\gamma_N]\in L|\; 
             s(\gamma_1(\zeta),\dots,\gamma_N(\zeta))=0\right\}.
\]
In general, it is the wrong space. 
In other words, 
its finite dimensional strata do not have the expected dimensions. 
Instead we regard it as a superspace. 
The coefficient of $\zeta^\nu$ of 
$s(\gamma_1(\zeta),\dots,\gamma_N(\zeta))$ takes values in 
the bundle $\bigoplus_{i=1}^l\mV_{i,\nu}$. 
Therefore, the superspace (\ref{eq:superspaceoverstablemanifold})
can be considered to approximate a stable manifold in $\widetilde{LY}$
when $d_1=d_2$.  

The module $A^*_{T^2}$ has a $D_{\rm poly}$-module structure, 
where $D_{\rm poly}$ is the polynomial part of $D$, i.e. 
$D_{\rm poly}\defequal \C[\Qq,\Qp,\hbar]$. 
The covering transformation $Q^a$ induces maps 
${Q^a}\colon (L,\mW_{d_1}^{d_2})\rightarrow (L,\mW_{d_1-q^b}^{d_2-q^b})$ 
and 
${Q^a}^*\colon H^*_{T^2}(L/\mW_{d_1-q^a}^{d_2-q^a})\rightarrow 
 H^*_{T^2}(L/\mW_{d_1}^{d_2})$ 
which is compatible with the inductive limit.  
Thus we have a pull back ${Q^a}^*\colon A^*_{T^2}\rightarrow A^*_{T^2}$
which raises the degree by $2\pair{c_1(\XSig/\mV)}{q^a}=\deg q^a$.  
By Proposition \ref{prop:commutationrelation}, 
$D_{\rm poly}$ acts on $A^*_{T^2}$ by $\Qq^a\mapsto {Q^a}^*$ 
and $\Qp_a\mapsto P_a\cup$.  
We define the $T^2$-equivariant 
semi-infinite cohomology $H^{\infty/2}_{T^2}$ 
as a submodule of $A^*_{T^2}$ by 
\[H^{\infty/2}_{T^2}\defequal \sum_{d\in H_2(\XSig,\Z)}
        H^*_{T^2}(L/\mW_{d}^{d}) \subset A^*_{T^2}. 
\]
It is a sub $D_{\rm poly}$-module of $A^*_{T^2}$. 
Define the Floer fundamental cycle $\Delta$ as the image 
in $H^{\infty/2}_{T^2}$ of 
the class $1$ in $H^*_{T^2}(L/\mW_{0}^{0})$. 
This class $\Delta$ is a fundamental class of 
the stable manifold $(L_{0}^\infty,\oplus_{i=1}^l
\oplus_{\nu \ge 0}\mV_{i,\nu})$. 

Semi-infinite homology can be defined by replacing  
$\hbar$ with $-\hbar$. 
Define 
\[\ov{\mW}_{d_1}^{d_2}\defequal 
  \bigoplus_{i=1}^l\bigoplus_{\nu=-\infty}^{\pair{v_i}{d_1}}
  \mathcal{V}_{i,\nu} \oplus 
  \bigoplus_{i=1}^N\bigoplus_{\nu=\pair{u_i}{d_2}+1}^{\infty} 
  \mathcal{L}_{i,\nu}.
\] 
Then we have an inductive system 
$\{H^*_{T^2}(L/\ov{\mW}_{d_1}^{d_2})\}_{d_1,d_2}$ and its limit
\[A_*^{T^2}\defequal \underset{d_1,d_2}{\varinjlim}
     H^*_{T^2}(L/\ov{\mW}_{d_1}^{d_2})
     \left[\textstyle -2\sum_{i=1}^N\pair{u_i}{d_2}+
           2\sum_{i=1}^l\pair{v_i}{d_1}
     \right].
\]
Then $A_*$ becomes a graded module because of the shift of grading. 
The covering transformation $Q^a$ gives a push-forward 
$Q^a_*\colon H_{T^2}^*(L/\ov{\mW}_{d_1}^{d_2})\rightarrow 
H_{T^2}^*(L/\ov{\mW}_{d_1-q^a}^{d_2-q^a})$. 
Therefore, $A_*^{T^2}$ 
has an action of $\ov{D}_{\rm poly}$, the polynomial part 
of $\ov{D}$, by $\Qq^a \mapsto Q^a_*$ and $\Qp_a\mapsto P_a\cup$.  
The operation $Q^a_*$ also raises the degree by $\deg q^a$. 
We define the $T^2$-equivariant semi-infinite homology 
$H_{\infty/2}^{T^2}$ as 
\[H_{\infty/2}^{T^2}\defequal \sum_{d} H^*_{T^2}(L/\ov{\mW}_{d}^{d}) 
           \subset A_*^{T^2}.
\]
This is preserved by the action of $\ov{D}_{\rm poly}$. 
Denote by $\ov{\Delta}$ 
the image in $H_{\infty/2}^{T^2}$ of the class $1$ in $H^*_{T^2}(L/\mW_0^0)$.  
We have a bar isomorphism $\barop$ between $A^*_{T^2}$ and $A_*^{T^2}$ 
defined by $\hbar\mapsto -\hbar$. 
More precisely, 
on the submodule $H^*_{T^2}(L/\mW_{d_1}^{d_2})$, 
it is defined by 
\[\barop\colon H^*_{T^2}(L/\mW_{d_1}^{d_2})
   \longrightarrow H^*_{T^2}(L/\ov{\mW}_{-d_2}^{-d_1}),\quad 
  f(P,\hbar)\mapsto f(P,-\hbar).  
\]
The bar isomorphism maps $H^{\infty/2}_{T^2}$ to $H_{\infty/2}^{T^2}$, 
i.e.  
$\barop\colon H^{\infty/2}_{T^2}\overset{\cong}
{\rightarrow}H_{\infty/2}^{T^2}$.

Next we define a pairing between $H^{\infty/2}_{T^2}$ 
and $H_{\infty/2}^{T^2}$. 
Let $\alpha$ and $\beta$ be elements of $H^*_{T^2}(L/\ov{\mW}_{d_2}^{d_2})$ 
and $H^*_{T^2}(L/\mW_{d_1}^{d_1})$ respectively. 
If $\pair{u_i}{d_2-d_1}\ge 0$ and $\pair{v_i}{d_2-d_1}\ge 0$ 
for all $i$, we define 
\[\int_{L}\alpha\cup\beta\defequal 
\int_{L_{d_1}^{\infty}\bigcap L_{-\infty}^{d_2}}\alpha\cup\beta
\cup\prod_{i=1}^l\prod_{\nu=\pair{v_i}{d_1}}^{\pair{v_i}{d_2}}V_{i,\nu}.
\]
This takes values in $\C[\hbar,\lambda]$. 
In general, however, it is not well-defined. 

\noindent
(Case 1): When there exists $i_0$ such that $\pair{v_{i_0}}{d_2-d_1}<0$.  
Because $v_{i_0}$ is nef, $d_2-d_1$ is not in the semigroup $\Lambda$.  
By Lemma \ref{lem:Moricone}, we have $\sigma_{I_0}\notin \Sigma$ 
for $I_0\defequal \{i|\;\pair{u_i}{d_2-d_1}<0\}$. 
Therefore the intersection $L_{d_1}^{\infty}\cap L_{-\infty}^{d_2}$ 
is empty and we set $\int_{L}\alpha\cup \beta=0$. 

\noindent
(Case 2): When $\pair{v_i}{d_2-d_1}\ge 0$ holds for all $i$. 
We choose $d'_2$ such that $\pair{d'_2-d_1}{u_i}\ge 0$ and 
$\pair{d'_2-d_2}{u_i}\ge 0$ hold for all $i$. 
Then we define 
\[\int_{L}\alpha\cup\beta\defequal 
\int_{L_{d_1}^{\infty}\bigcap L_{-\infty}^{d'_2}}\alpha\cup\beta
\cup\prod_{i=1}^l\prod_{\nu=\pair{v_i}{d_1}}^{\pair{v_i}{d_2}}V_{i,\nu}
\cup\prod_{i=1}^N\prod_{\nu=\pair{u_i}{d_2}+1}^{\pair{u_i}{d'_2}}U_{i,\nu}.
\]
This does not depend on the choice of $d'_2$. 

We also define a $q$-deformed pairing $(\cdot,\cdot)$ between 
semi-infinite homology and cohomology as follows: 
\begin{equation}
\label{eq:definitionofpairinginFloertheory}
(\alpha,\beta)\defequal \sum_{d\in H_2(\XSig,\Z)}
          q^d\int_{L}\alpha\cup {Q^{-d}}^*(\beta),\quad 
  \text{for }\alpha\in H_{\infty/2}^{T^2}, \beta\in H^{\infty/2}_{T^2}, 
\end{equation}
where 
$Q^d\defequal
(Q^1)^{\pair{p_1}{d}}(Q^2)^{\pair{p_2}{d}}\cdots(Q^r)^{\pair{p_r}{d}}$. 
It decreases the degree by $\dim_{\R} \XSig=2n$ and  
takes values in $\C[\hbar,\lambda][\![q,q^{-1}]\!]$.

Define $\Delta_d\defequal (Q^d)^*\Delta$ for $d$ in $H_2(\XSig,\Z)$. 
Let $i_d\colon (\XSig)_d\hookrightarrow L$ be the inclusion. 
Hereafter we consider $\Delta$ and $\Delta_d$ 
to be expanded as infinite products, i.e. 
\[\Delta_d=\prod_{i=1}^N\prod_{\nu<\pair{u_i}{d}}U_{i,\nu}\cup
           \prod_{i=1}^l\prod_{\nu\ge\pair{v_i}{d}}V_{i,\nu}.
\]
We define a localization map 
$i_d^*(\cdot/\Delta_d)\colon A^*_{T^2}\rightarrow 
 H^*(\XSig)\otimes\C(\lambda,\hbar)$ 
by using Lemma \ref{lem:UVPrelation} as 
\begin{align}
\label{eq:localizationmap} 
\nonumber
i_d\left(\frac{\alpha}{\Delta_d}\right)  
   & \defequal i_d^*\left(\alpha\cup\prod_{i=1}^N
                \frac{\prod_{\nu<\pair{u_i}{d_1}}U_{i,\nu}}
                     {\prod_{\nu<\pair{u_i}{d}}U_{i,\nu}}\cup 
                           \prod_{i=1}^l
                \frac{\prod_{\nu\ge\pair{v_i}{d_2}}V_{i,\nu}}
                     {\prod_{\nu\ge\pair{v_i}{d}}V_{i,\nu}}\right) \\
   & \defequal i_d^*(\alpha)\cup\prod_{i=1}^N
                \frac{\prod_{\nu>0}(u_i+(\pair{u_i}{d-d_1}+\nu)\hbar)}
                     {\prod_{\nu>0}(u_i+\nu\hbar)}\cup 
          \prod_{i=1}^l
          \frac{\prod_{\nu\ge0}(v_i+(\pair{v_i}{d-d_2}-\nu)\hbar-\lambda)}
                     {\prod_{\nu\ge0}(v_i-\nu\hbar-\lambda)} 
\end{align}  
for $\alpha\in H^*_{T^2}(L/\mW_{d_1}^{d_2})$. 
We can see that the right hand side gives an element in 
$H^*(\XSig)\otimes\C(\lambda,\hbar)$. 
The denominator $i_d^*(\Delta_d)$ 
is considered to be the Euler class of 
the negative normal bundle of $(\XSig)_d$. 
Therefore, this map defines a partial integration over an open stratum of  
the unstable manifold 
$(L_{-\infty}^d,\oplus_{i=1}^l \oplus_{\nu\le\pair{u_i}{d}} \mV_{i,\nu})$
of $(\XSig)_d$. 
This localization also appears in \cite{vlassopoulos}.  
We next describe the restriction of the localization map to 
$H^{\infty/2}_{T^2}$. 

\begin{lemma}
On the submodule $H^{\infty/2}_{T^2}$, 
$i_d^*(\cdot/{\Delta_d})$ takes values 
in $H^*(\XSig)\otimes\C[\lambda,\hbar,\hbar^{-1}]$. 
Moreover, if $d\notin\Lambda$, we have 
$i_d^*(\Delta/\Delta_d)=0$. 
\end{lemma}
\begin{proof}
Consider the case where $d_1=d_2$ in the equation 
(\ref{eq:localizationmap}), 
and set $d'=d-d_1=d-d_2$.  
It suffices to consider the case where 
there exists $i_0$ such that $\pair{v_{i_0}}{d'}<0$ holds. 
In this case, $d'$ is not in the semigroup $\Lambda$ 
because $v_{i_0}$ is nef.  
Therefore by Lemma \ref{lem:Moricone}, we have
$\sigma_{I_0}\notin\Sigma$ for 
$I_0\defequal \{i|\;\pair{u_i}{d'}<0\}$. 
Hence by (\ref{eq:relationoftoric}), we have 
$\prod_{i\in I_0}u_i=0$. 
Because the first term contains this factor, 
$i_d^*(\alpha/\Delta_d)$ vanishes. 
The last statement follows from the same argument. 
\end{proof}

By using the lemma, we define a map 
$\Xi\colon H^{\infty/2}_{T^2}\rightarrow H^*(\XSig)\otimes 
\C[\hbar,\hbar^{-1},\lambda][\![q,q^{-1}]$ by 
\begin{equation}
\label{eq:definitionofXi}
\Xi(\alpha)\defequal \sum_{d\in H_2(\XSig,\Z)}
             q^di_d^*\left(\frac{\alpha}{\Delta_d}\right).
\end{equation}
\begin{proposition} 
The map $\Xi$ is a $\C[q^1,\dots,q^r]$-linear 
injection and preserves the degree. 
Moreover, restricted on $\C[Q^*,P,\hbar,\lambda]\Delta$, 
it takes values in $H^*(\XSig)\otimes \C[\hbar,\hbar^{-1},\lambda][\![q]\!]$. 
\end{proposition}
\begin{proof}
It is clear that $\Xi$ preserves the degree. 
First we have 
\[i_d^*\left(\frac{{Q^a}^*(\alpha)}{\Delta_d}\right)
       =i_d^*{Q^a}^*\left(\frac{\alpha}{\Delta_{d-q^a}}\right)
       =i_{d-q^a}^*\left(\frac{\alpha}{\Delta_{d-q^a}}\right)
\]
because $Q^a\circ i_d=i_{d-q^a}$. 
Thus, $\Xi$ is $\C[q^1,\dots,q^r]$-linear. 
Set $\alpha=f(P,\hbar,\lambda)\Delta$ for some polynomial $f$.  
If $\pair{p_a}{d}<0$ holds for some $a$, we have 
$i_d^*(\Delta/\Delta_d)=0$ because $p_a$ is nef and 
$d\notin \Lambda$. 
Therefore $\Xi(\alpha)$ has no negative powers of $q^a$. 
From this we see that $\Xi(\C[Q^*,P,\hbar,\lambda]\Delta)\subset 
H^*(\XSig)\otimes \C[\hbar,\hbar^{-1},\lambda][\![q]\!]$.  
Finally we show the injectivity of $\Xi$. 
Suppose we have $\Xi(\alpha)=0$ for 
$\alpha=f(P,\hbar,\lambda)\in H^*_{T^2}(L/\mW_{d_1}^{d_2})$. 
Take $d'$ such that $\pair{u_i}{d'-d_j}\ge 0$ holds for all $i$, $j$. 
In this case we also have $\pair{v_i}{d'-d_j}\ge0$ 
because $d'-d_j$ is in $\Lambda$ by Lemma \ref{lem:Moricone}.  
Then we have 
\[0=i_{d'}^*\left(\frac{\alpha}{\Delta_{d'}}\right)
   =f(p_a+\pair{p_a}{d'}\hbar,\hbar,\lambda)\cup 
    \frac{\prod_{i=1}^l\prod_{\nu=1}^{\pair{v_i}{d'-d_2}}
          (v_i+\nu\hbar-\lambda)}
         {\prod_{i=1}^N\prod_{\nu=1}^{\pair{u_i}{d'-d_1}}
          (u_i+\nu\hbar)}.
\]
Therefore, we have $f(p_a+\pair{p_a}{d'}\hbar,\hbar,\lambda)=0$ 
in $H^*(X)\otimes \C[\hbar,\lambda]$. 
It holds for all $d'$ satisfying $\pair{u_i}{d'-d_j}\ge 0$, 
hence we have $f=0$ as a polynomial. 
\end{proof}

Lemma \ref{lem:UVPrelation} immediately establishes 
the following differential equation 
which is the same as the one 
in Lemma \ref{lem:differentialequationforXi}. 
\begin{proposition}
\label{prop:differentialequationforXiinFloertheory}
\[
\Xi(P_a\cup\alpha)=(\hbar\parti{a}+p_a)\Xi(\alpha).
\]
\end{proposition}

In semi-infinite homology, 
we define $\ov{\Delta}_d\defequal Q^d_*(\ov{\Delta})$. 
Similarly, for an element $\alpha$ in $H_{\infty/2}^{T^2}$, 
we have a well-defined element $i_{-d}^*(\alpha/\ov{\Delta}_d)$ 
in $H^*(\XSig)\otimes \C[\lambda,\hbar,\hbar^{-1}]$.  
Therefore we can define 
$\ovXi\colon H_{\infty/2}^{T^2}\rightarrow H^*(\XSig)\otimes
\C[\lambda,\hbar,\hbar^{-1}][\![q,q^{-1}]$ by  
\[\ovXi(\alpha)\defequal \sum_{d\in H_2(\XSig,\Z)}q^d i_{-d}^*
         \left(\frac{\alpha}{\ov{\Delta}_d}\right). 
\]
This is related to $\Xi$ by the formula 
$\ov{\Xi}=\barop\circ\Xi\circ\barop$, 
where $\barop$ acts on 
$H^*(\XSig)\otimes \C[\hbar,\hbar^{-1},\lambda][\![q]\!]$ 
by $\hbar\mapsto -\hbar$. 
The map $\ovXi$ is also a $\C[q^1,\dots,q^r]$-linear injection and 
preserves the degree. 
Moreover, restricted on $\C[Q_*,P,\hbar,\lambda]\Delta$, it takes 
values in $H^*(\XSig)\otimes\C[\hbar,\hbar^{-1},\lambda][\![q]\!]$.  

We calculate the pairing $\int_{L}$ by using the 
localization theorem in equivariant cohomology, 
see e.g.  \cite{atiyah-bott}. 
For $d$ in $H_2(\XSig,\Z)$, we define a map  
$i_d^*(\cdot\cup\cdot)/E_d$ as follows: 
\begin{gather*}
\frac{i_d^*(\alpha\cup\beta)}{E_d}
  \defequal 
  i_d^*\left(\frac{\alpha}{\ov{\Delta}_{-d}}\right)
  i_d^*\left(\frac{\beta}{\Delta_d}\right)\cup \prod_{i=1}^l(v_i-\lambda), 
\quad 
\text{for }\alpha\in H_{\infty/2}^{T^2}, \beta\in H^{\infty/2}_{T^2}.
\end{gather*}
Formally, $E_d$ can be written as 
$\prod_{i=1}^N\prod_{\nu\neq 0}(u_i-\nu\hbar)
 \prod_{i=1}^l\prod_{\nu}(v_i-\nu\hbar-\lambda)$. 
It is considered to be the Euler class of the normal bundle of $(\XSig)_d$. 
A direct application of the localization theorem shows 
\begin{proposition}[localization]
\[
\int_{L}\alpha\cup\beta=\sum_{d\in H_2(\XSig,\Z)}\int_{\XSig}
         \frac{i_d^*(\alpha\cup\beta)}{E_d}.
\]
\end{proposition}
This is an infinite dimensional version of the localization theorem. 
Note that each term on the right hand side 
is in the localized ring $\C[\lambda,\hbar,\hbar^{-1}]$, 
but the total sum is in $\C[\lambda,\hbar]$.

\begin{proposition}
\label{prop:pairingofFloercohomology}
For $\alpha\in H_{\infty/2}^{T^2}(L)$ 
and $\beta\in H^{\infty/2}_{T^2}(L)$, 
we have 
\[
(\alpha,\beta)=\int_{\XSig}\ov{\Xi}(\alpha)\cup \Xi(\beta)
\cup Euler_{S^1}(\mV). 
\]
\end{proposition} 
\begin{proof}
This is a direct consequence of the localization formula. 
\begin{eqnarray*}
\allowdisplaybreaks 
(\alpha,\beta)&=&\sum_{d\in H_2(\XSig,\Z)}q^d
        \sum_{d'\in H_2(\XSig,\Z)}\int_{\XSig}
        \frac{i_{d'}^*(\alpha\cup{Q^{-d}}^*\beta)}{E_{d'}} \\
  &=& \sum_{d\in H_2(\XSig,\Z)}q^d
      \sum_{d'\in H_2(\XSig,\Z)}\int_{\XSig}
         i_{d'}^*\left(\frac{\alpha}{\ov{\Delta}_{-d'}}\right)
      \cup i_{d'}^*\left(\frac{{Q^{-d}}^*(\beta)}{\Delta_{d'}}\right)
      \cup Euler_{S^1}(\mV) \\  
  &=& \sum_{d,d'}\int_{\XSig}
      q^{-d'}i_{d'}^*\left(\frac{\alpha}{\ov{\Delta}_{-d'}}\right)
      \cup q^{d+d'}i_{d+d'}^*\left(\frac{\beta}{\Delta_{d+d'}}\right)
      \cup Euler_{S^1}(\mV) \\
  &=& \int_{\XSig}\ov{\Xi}(\alpha)\cup \Xi(\beta)\cup Euler_{S^1}(\mV).
\end{eqnarray*} 
\end{proof}
The proposition corresponds to the 
definition of pairing in the abstract quantum $D$-module. 

Finally we define equivariant Floer cohomology and homology. 
Consider a $D_{\rm poly}$ module 
$H^{\infty/2}_{T^2,+}\defequal 
\C[Q^*,P,\hbar,\lambda]\Delta\subset H^{\infty/2}_{T^2}$ 
generated by $\Delta$. 
It has a $Q$-adic topology which is Hausdorff 
because of the injectivity of $\Xi$. 
Define equivariant Floer cohomology $FH^*_{T^2}$ by
\[FH_{T^2}^*\defequal \widehat{H}^{\infty/2}_{T^2,+}
\defequal \underset{n}{\varprojlim} 
H^{\infty/2}_{T^2,+}/\mathfrak{m}^nH^{\infty/2}_{T^2,+}, 
\]
where $\mathfrak{m}$ is the ideal of $D_{\rm poly}$ generated by ${Q^a}^*$.  
Then $FH^*_{T^2}$ becomes a $D$-module (not just a $D_{\rm poly}$-module).  
The map $\Xi$ can be extended to this completion and 
defines an embedding 
$\Xi\colon FH^*_{T^2}\hookrightarrow 
 H^*(\XSig)\otimes\C[\hbar,\hbar^{-1},\lambda][\![q]\!]$. 
Similarly, we define the equivariant Floer homology $FH_*^{T^2}$. 
We have a bar isomorphism 
$\barop\colon FH^*_{T^2}\overset{\cong}{\rightarrow}FH_*^{T^2}$.    
$FH_*^{T^2}$ is a $\ov{D}$-module and we have an embedding 
$\ovXi\colon  FH_*^{T^2}\hookrightarrow 
 H^*(\XSig)\otimes\C[\hbar,\hbar^{-1},\lambda][\![q]\!]$. 
By Proposition \ref{prop:pairingofFloercohomology}, 
we extend the pairing $(\cdot,\cdot)$ 
between Floer homology and cohomology.   
\[(\cdot,\cdot)\colon FH^{T^2}_*\times FH^*_{T^2}
\longrightarrow \C[\hbar,\lambda][\![q]\!].
\]

\begin{remark}
So far, we have studied $T^2$-equivariant Floer cohomology. 
In the case of a toric variety itself,  
there is no $S^1$ acting on the fiber of $\mV$. 
Hence, it is natural to consider 
$S^1$-equivariant Floer cohomology $FH^*_{S^1}$. 
(There remains the $S^1$ action rotating loops).  
For the superspace, 
the fiberwise $S^1$ action on $\mV$ is introduced in order to ensure that 
the pairing becomes non-degenerate. 
After we take a frame $\Phi$ of $T^2$-equivariant Floer cohomology 
as in the next subsection, 
we can take $\lambda$ to be zero and 
obtain the $S^1$-equivariant version. 
In Section \ref{sect:abstractquantumDmodule}, 
we considered only the case where $\lambda = 0$, however, 
it easy to develop the theory including $\lambda$. 
\end{remark}

\subsection{Quantum D-modules and Equivariant Floer Theory} 

We shall show that the equivariant 
Floer cohomology is an abstract quantum $D$-module 
and coincides with the quantum $D$-module arising from the quantum 
cohomology of the superspace $(\XSig,\mV)$.  
First, we must define its frame. 
By the existence of the embedding $\Xi$, 
the zero-fiber of $FH^*_{T^2}$ is isomorphic to 
$H^*(\XSig)\otimes\C[\hbar,\lambda]$ 
through the localization map $i_0^*(\cdot/\Delta)$. 
We regard this isomorphism as a fixed choice of a frame of 
the zero-fiber. 
Choose a basis $\{T_0(p),\dots,T_s(p)\}$ of $H^*(\XSig)$,
where $T_i$ is a polynomial in $p_1,\dots,p_r$ and $T_0(p)=1$, 
$T_a(p)=p_a$ for $1\le a\le r$. 
Then we define a frame $\Phi$ as the  
$\C[\hbar,\lambda][\![q]\!]$-linear map:  
\[\Phi\colon H^*(\XSig)\otimes\C[\hbar,\lambda][\![q]\!] 
 \rightarrow FH^*_{T^2}, \quad 
 T_a(p)\mapsto T_a(P)\Delta. 
\] 
This clearly preserves the degree. 
Correspondingly, we define $\ovPhi$ by 
\[\ovPhi\colon H^*(\XSig)\otimes\C[\hbar,\lambda][\![q]\!] 
 \rightarrow FH_*^{T^2}, \quad 
 T_a(p)\mapsto T_a(P)\ov{\Delta}. 
\] 
First we claim that $\Phi$ defines an isomorphism 
over  $\C[\hbar,\lambda][\![q]\!]$. 
Composing $\Phi$ with $\Xi$, we have 
\begin{align}
\nonumber
\Xi\circ \Phi\colon &H^*(\XSig)\otimes 
   \C[\hbar,\lambda][\![q]\!]\rightarrow 
  H^*(\XSig)\otimes \C[\hbar,\hbar^{-1},\lambda][\![q]\!], 
\\ 
\label{eq:XiPhi}
  & T_a(p) \mapsto  T_a(p)+\sum_{d\neq 0}
      i_d^*
      \left(\frac{T_a(P)\Delta}{\Delta_d}\right)q^d.
\end{align}
From this we see that $\Xi\circ\Phi$ is injective. 
Therefore, $\Phi$ is also injective. 
For the surjectivity, 
it suffices to show that $f(P)\cdot\Delta$ is 
in the image of $\Phi$ for any polynomial $f$.   
By (\ref{eq:XiPhi}), 
we can easily deduce that there exists an element $x_0$ in 
$H^*(\XSig)\otimes \C[\hbar,\hbar^{-1},\lambda][\![q]\!]$ such that 
$\Xi\circ\Phi(x_0)=\Xi(f(P)\Delta)$ holds.  
We would like to show that $x_0$ is an element of 
$H^*(\XSig)[\hbar,\lambda][\![q]\!]$.  
Consider the pairing $(\cdot,\cdot)_{\Phi}$ defined by 
\begin{eqnarray*}
(\alpha,\beta)_{\Phi}&\defequal& 
  (\ovPhi(\alpha),\Phi(\beta)) \\
 &=& \int_{\XSig}\alpha\cup\beta\cup Euler_{S^1}(\mV)
     + \text{higher order terms in $q$}.
\end{eqnarray*}
If we invert the variable $\lambda$, 
this becomes a perfect pairing on 
$H^*(\XSig)\otimes\C[\hbar,\lambda,\lambda^{-1}][\![q]\!]$
which takes values in $\C[\hbar,\lambda,\lambda^{-1}][\![q]\!]$. 
On the other hand, we have for any $a$, 
\[(T_a(p),x_0)_{\Phi}
  =(\ov{\Phi}(T_a(p)),f(P)\Delta)\in \C[\lambda,\hbar][\![q]\!]. \] 
Therefore, $x_0$ is contained in 
$H^*(\XSig)\otimes\C[\hbar,\lambda,\lambda^{-1}][\![q]\!]$. 
Thus we have $x_0\in H^*(\XSig)\otimes \C[\hbar,\lambda][\![q]\!]$ 
and the surjectivity of $\Phi$ is proved.  

It is easy to check that $\Phi$ induces the frame 
$\Phi_0=i_0(\cdot/\Delta)$ of the zero fiber and 
$\Phi_0$ becomes an isomorphism of $\C[p,\hbar]$-modules. 
By definition, we have $\Phi(p_a)=P_a\Phi(1)$. 
Therefore $\Phi$ and the generalized coordinates $(Q^*,P)$ are compatible. 
This frame defines a connection 
$\nabla=d+\frac{1}{\hbar}\sum\Omega_a dq^a/q^a$ on the 
$H^*(\XSig)$ bundle. 
It satisfies $\lim_{q\to 0}\Omega_a=(p_a\cup)$. 

Next we claim that the $\End(V)$-valued function  
$S$ in Proposition \ref{prop:existenceofS} 
can be written as $S\defequal(\Xi\circ\Phi)^{-1}$. 
The initial condition $S(q=0,\hbar)=I$ is clear. 
Proposition \ref{prop:differentialequationforXiinFloertheory} 
shows that $S$ satisfies the differential equation: 
\[ \hbar\nabla^1_aS=\hbar\nabla_aS-S\circ(p_a\cup)=0.
\] 
The equation 
$\nabla^1_{\Eul}S=[\adE,S]=0$ holds because 
$S$ preserves degrees. 
Consequently $S$ is the flat section in 
Proposition \ref{prop:existenceofS}. 
In particular, $\Xi$ in this section coincides with 
that in Section \ref{sect:abstractquantumDmodule}. 
Therefore, by Proposition \ref{prop:pairingofFloercohomology}, 
the pairing as an abstract quantum $D$-module is identical with 
the pairing $(\cdot,\cdot)$ 
between equivariant Floer homology and cohomology. 
The pairing $(\cdot,\cdot)$ in the Floer theory 
is polynomial in $\hbar$. 
This fact together with Proposition \ref{prop:polynomiality}
and Remark \ref{rem:pairing} shows that 
$(\cdot,\cdot)$ can be identified with 
the Poincar\'{e} pairing $\pairV{\cdot}{\cdot}$ 
of the superspace $(\XSig,\mV)$ 
after we take the canonical frame.   

The $J$-function of $FH^*_{T^2}$ is calculated as follows. 
\begin{eqnarray*}
J_{\Phi,(Q^*,P)} &=& e^{p\log q/\hbar}S^{-1}(1)
 =e^{p\log q/\hbar}\sum_{d;\pair{p_a}{d}\ge 0,\forall a}
  i_d^*\left(\frac{\Delta}{\Delta_d}\right)q^d \\
&=& e^{p\log q/\hbar} \sum_{d;\pair{p_a}{d}\ge 0,\forall a}
             q^d\prod_{i=1}^N 
       \frac{\prod_{\nu=\pair{u_i}{d}+1}^{\infty}(u_i+\nu\hbar)}
            {\prod_{\nu=1}^{\infty}(u_i+\nu\hbar)}\cup 
   \prod_{i=1}^l
      \frac{\prod_{\nu=-\infty}^{\pair{v_i}{d}}(v_i+\nu\hbar-\lambda)}
           {\prod_{\nu=-\infty}^{0}(v_i+\nu\hbar-\lambda)}.  
\end{eqnarray*}
By Givental \cite{givental-mirrorthm-toric}, 
when $c_1(\XSig/\mV)$ is nef, 
the above function coincides with the $J$-function
of the superspace $(\XSig,\mV)$ after a suitable mirror transformation. 
The necessary coordinate change can be done uniquely as in Proposition 
\ref{prop:giventalmirrortrans}. 
Because abstract quantum $D$-modules are determined 
by $J$-functions (Theorem \ref{thm:abstractreconstructionbyJ}), 
we conclude that our equivariant Floer cohomology 
$FH^*_{T^2}$ is isomorphic to
the quantum $D$-module of the superspace $(\XSig,\mV)$ 
as an abstract quantum $D$-module. 
We remark that there exist natural generalized coordinates 
$({Q^1}^*,\dots,{Q^r}^*,P_1,\dots,P_r)$ in equivariant Floer theory. 
These coordinates are not necessarily affine, 
and we need a mirror transformation. 
Even in that case, however, we need not change the normalization, i.e. 
the normalization vector of equivariant Floer theory is identical with  
that of the quantum $D$-module.

\begin{theorem}
Let $\XSig$ be a toric manifold
and $\mV=\bigoplus_{i=1}^l\mV_i$ be a sum of nef line bundles 
over $\XSig$. 
The equivariant Floer cohomology $FH^*_{T^2}$ is an 
abstract quantum $D$-module 
endowed with a choice of a frame $\Phi_0$ of the zero-fiber 
and generalized coordinates $({Q^1}^*,\dots,{Q^r}^*,P_1,\dots,P_r)$. 
If $c_1(\XSig/\mV)\defequal c_1(\XSig)-c_1(\mV)$ is nef, 
it is isomorphic to the quantum $D$-module of the 
superspace $(\XSig,\mV)$ 
as an abstract quantum $D$-module with a fixed choice of $\Phi_0$ 
and a normalization vector 
$v_0=(\partial/\partial q^1)_0+\dots+(\partial/\partial q^r)_0$.  
Moreover, 

$(1)$ The map $\Xi$ defined in $(\ref{eq:definitionofXi})$
coincides with the $\Xi$ as an abstract quantum $D$-module. 

$(2)$ The pairing $(\cdot,\cdot)$ 
between equivariant Floer homology and cohomology 
defined in $(\ref{eq:definitionofpairinginFloertheory})$ 
coincides with that as an abstract quantum $D$-module. 
It can be identified with the Poincar\'{e} pairing 
$\pairV{\cdot}{\cdot}$ after we take the canonical frame.  
\end{theorem} 

\begin{remark}
From the above theorem and Corollary \ref{cor:JfunctiondeterminesGW}, 
we can calculate all genus zero Gromov-Witten invariants of $\XSig$ 
by using the equivariant Floer cohomology. 
\end{remark}


\section{Examples} 

We illustrate the general theory with 
two examples concerning a Hirzeburch surface. 
We obtain the quantum multiplication tables by 
using the equivariant Floer theory constructed in 
Section \ref{sect:equivariantFloertheory}. 

Consider the Hirzeburch surface  
$\F_1\defequal \Proj(\mathcal{O}(1)\oplus \mathcal{O})$, 
where $\mathcal{O}$ is the structure sheaf of $\Proj^1$. 
It is described by the fan $\Sigma$ in $\R^2$ 
whose one skeleton $\Sigma^{(1)}$ consists of 
four vectors 
\[x_1=(1,0),\quad x_2=(-1,-1),\quad x_3=(0,1),\quad 
  x_4=(0,-1). 
\]
Let $u_i$ be the cohomology class corresponding to $x_i$. 
We have the following relations.  
\[u_1=u_2, \quad u_3=u_4+u_2,\quad u_1u_2=0,\quad u_3u_4=0. 
\]
We can see that 
$u_1$ is the class of the fiber of $\F_1\rightarrow \Proj^1$ and that 
$u_3$ and $u_4$ are classes of the zero-section and the
$\infty$-section respectively.  
We take $p_1\defequal u_1$ and $p_2\defequal u_3$ as a basis 
of $H^2(\F_1)$. 
The K\"{a}hler cone of $H^2(\F_1)$ is generated by 
$p_1$ and $p_2$. 
The set $\{1,p_1,p_2,p_1p_2\}$ forms a linear basis 
of the total cohomology ring $H^*(\F_1)$.  
The first Chern class of $\F_1$ is $c_1(\F_1)=p_1+2p_2$.  
The classes $p_1$, $p_2$ may be lifted to 
the equivariant classes $P_1$, $P_2$ in  
Givental's model for $\widetilde{L\F_1}$. 
Let $Q^1$, $Q^2$ denote the covering transformations dual to 
$p_1$, $p_2$ and  
$q^1$, $q^2$ denote the corresponding variables. 
We have $\deg q^1=2$ and $\deg q^2=4$. 
The Floer fundamental cycle is  
\[\Delta=\prod_{\nu<0}(P_1-\nu\hbar)^2\prod_{\nu<0}(P_2-\nu\hbar)
         \prod_{\nu<0}(P_2-P_1-\nu\hbar).
\]
From this and the formula ${Q^a}^*(P_b)=P_b-\delta^a_b\hbar$, 
we can derive the Picard-Fuchs equations. 
\[(P_2-P_1+\hbar){Q^1}^*(\Delta)=P_1^2\Delta, \quad 
  {Q^2}^*(\Delta)=P_2(P_2-P_1)\Delta.
\]
The equivariant Floer cohomology $FH_{S^1}^*$ 
is generated by $\Delta$ over the Heisenberg algebra 
$\C[\hbar][\![Q^*]\!][P]$ 
with the above relations.  
We take a frame $\Phi$ of $FH_{S^1}^*$ 
as 
\[ \Phi(1)\defequal \Delta, \quad 
   \Phi(p_1)\defequal P_1\Delta, \quad 
   \Phi(p_2)\defequal P_2\Delta, \quad 
   \Phi(p_1p_2)\defequal P_1P_2\Delta.
\]
By using the Picard-Fuchs equations, we easily obtain the connection
matrices $\Omega_1$, $\Omega_2$ 
defined by $\Phi(\Omega_a(\alpha))=P_a\Phi(\alpha)$. 
\[\Omega_1=
  \begin{pmatrix}
  0 & 0    & 0 & q^1q^2 \\
  1 & -q^1 & 0 & 0      \\
  0 & q^1  & 0 & 0      \\
  0 & 0    & 1 & 0      
  \end{pmatrix}, \quad 
  \Omega_2=
  \begin{pmatrix}
  0 & 0 & q^2 & q^1q^2 \\
  0 & 0 & 0   & q^2    \\
  1 & 0 & 0   & 0      \\
  0 & 1 & 1   & 0 
  \end{pmatrix}. 
\]
Note that $\Phi(q^a\alpha)={Q^a}^*(\Phi(\alpha))$ 
by definition. 
These connection matrices are already $\hbar$-independent, 
therefore the frame $\Phi$ defined above is canonical. 
Hence we can identify the above matrices with 
the quantum multiplications by $p_1$, $p_2$, i.e. $p_a*= \Omega_a$. 
They agree with Appendix 2 in \cite{guest-introd1}. 

Next we consider the case where we need a mirror transformation.  
Let $\mV$ be a semi-positive line bundle over $\F_1$ such that 
$c_1(\mV)=2p_2$. 
Then we have a convex superspace $(\F_1,\mV)$ whose 
first Chern class is $c_1(\F_1/\mV)=p_1$. 
In this case, we have $\deg q^1=2$ and $\deg q^2=0$. 
The Floer fundamental cycle $\Delta^{\mV}$ is 
written as 
\[\Delta^{\mV}=\prod_{\nu<0}(P_1-\nu\hbar)^2\prod_{\nu<0}(P_2-\nu\hbar)
         \prod_{\nu<0}(P_2-P_1-\nu\hbar)\prod_{\nu\ge 0}(2P_2-\nu\hbar). 
\]
From this, we derive the Picard-Fuchs equations as 
\begin{equation}
\label{eq:Picard-Fuchsofsuperspace}
(P_2-P_1+\hbar){Q^1}^*\Delta^{\mV}=P_1^2\Delta^{\mV}, \quad
  2P_2(2P_2-\hbar){Q^2}^*\Delta^{\mV}=P_2(P_2-P_1)\Delta^{\mV}.
\end{equation}
The equivariant Floer cohomology ${FH^*_{S^1}}^\mV$ is 
generated by $\Delta^{\mV}$ with the above relations. 
We define a frame $\Phi$ of ${FH^*_{S^1}}^\mV$ 
in the same way. 
Then, by using the Picard-Fuchs equations for $\Delta^{\mV}$, 
we obtain connection matrices $\Omega_1$, $\Omega_2$
as follows. 
\[\Omega_1=
  \begin{pmatrix}
  0 &  0 & 0 & 2\hbar^2\frac{xy}{1-4y} \\
  1 & -x & 0 & 0                       \\
  0 &  x & 0 & 6\hbar\frac{xy}{1-4y}   \\
  0 &  0 & 1 & 4\frac{xy}{1-4y}         
  \end{pmatrix}, \quad  
  \Omega_2=
  \begin{pmatrix}
  0 & 0 & 2\hbar^2\frac{y}{1-4y} & 2\hbar^2\frac{xy}{(1-4y)^2} \\
  0 & 0 & 0                          & 2\hbar^2\frac{y}{1-4y}        \\
  1 & 0 & 6\hbar\frac{y}{1-4y}   & 6\hbar\frac{xy}{(1-4y)^2}   \\
  0 & 1 & \frac{1}{1-4y}           & 6\hbar\frac{y}{1-4y}+
                                       4\frac{xy}{(1-4y)^2}  
  \end{pmatrix}.
\]
Here, we put $x=q^1$, $y=q^2$ for notational convenience. 
These connection matrices contain $\hbar$, therefore 
$\Phi$ is not canonical in this case. 
In order to obtain a canonical frame, we perform 
the Birkhoff factorization for $S_{\Phi}$. 
According to Section \ref{sect:equivariantFloertheory}, 
the inverse of $S_{\Phi}$ is given by 
the composition $\Xi\circ\Phi$.  
\begin{gather*}
S_{\Phi}^{-1}(\alpha)=\Xi(\Phi(\alpha))
     =\sum_{d_1,d_2\ge 0}
     x^{d_1} y^{d_2} i_{d_1,d_2}^*\left(
     \frac{\alpha\Delta^{\mV}}{({Q^1}^*)^{d_1}({Q^2}^*)^{d_2}\Delta^{\mV}}
     \right ),  
\end{gather*}
where $i^*_{d_1,d_2}(P_a)=p_a+d_a\hbar$. 
The $J$-function $J_{\Phi}$ is given by 
\begin{align*}
e^{-(p_1\log x+p_2\log y)/\hbar}J_{\Phi}&=
                S_{\Phi}^{-1}(1)= 
      \sum_{d_1,d_2\ge 0} 
      \frac{x^{d_1}y^{d_2}\prod_{\nu=1}^{2d_2}(2p_2+\nu\hbar)
            \prod_{\nu=d_2-d_1+1}^{\infty}(p_2-p_1+\nu\hbar)}
           {\prod_{\nu=1}^{d_1}(p_1+\nu\hbar)^2
            \prod_{\nu=1}^{d_2}(p_2+\nu\hbar)
            \prod_{\nu=1}^{\infty}(p_2-p_1+\nu\hbar)}  \\  
   &= I_0+I_1\frac{p_1}{\hbar}+I_2\frac{p_2}{\hbar}
             +I_3\frac{p_1p_2}{\hbar^2},  
\end{align*}
where $I_0$, $I_1$, $I_2$, $I_3$ are defined as follows:  
\begin{align*}
 I_0&\defequal\sum_{d_2\ge d_1}k_{d_1,d_2}
              \frac{x^{d_1}}{\hbar^{d_1}}y^{d_2},\quad  \\ 
 I_1& \defequal\sum_{d_2\ge d_1}k_{d_1,d_2} 
                (A_{d_2-d_1}-2A_{d_1})\frac{x^{d_1}}{\hbar^{d_1}}y^{d_2} 
    +\sum_{d_2<d_1}l_{d_1,d_2} 
              (-1)^{d_1-d_2}\frac{x^{d_1}}{\hbar^{d_1}}y^{d_2}, \\
 I_2&\defequal\sum_{d_2\ge d_1}k_{d_1,d_2}
          (2A_{2d_2}-A_{d_2}-A_{d_2-d_1})\frac{x^{d_1}}{\hbar^{d_1}}y^{d_2}
      -\sum_{d_2<d_1}l_{d_1,d_2} 
              (-1)^{d_1-d_2}\frac{x^{d_1}}{\hbar^{d_1}}y^{d_2}, \\
 I_3&\defequal\sum_{d_2\ge d_1}k_{d_1,d_2}
               \frac{x^{d_1}}{\hbar^{d_1}}y^{d_2}D_{d_1,d_2}
      +\sum_{d_2<d_1}l_{d_1,d_2}(2A_{d_1}-A_{d_1-d_2-1})
              (-1)^{d_1-d_2}\frac{x^{d_1}}{\hbar^{d_1}}y^{d_2}, \\  
k_{d_1,d_2}&\defequal \frac{(2d_2)!}{(d_1!)^2d_2!(d_2-d_1)!},
\quad 
l_{d_1,d_2} \defequal \frac{(2d_2)!(d_1-d_2-1)!}{(d_1!)^2d_2!}, \\
D_{d_1,d_2} &\defequal 4B_{2d_2}+B_{d_2}+C_{d_2}-B_{d_2-d_1}-C_{d_2-d_1}
   -4A_{d_1}A_{2d_2}-2A_{2d_2}A_{d_2}+2A_{d_1}A_{d_2}
     +2A_{d_1}A_{d_2-d_1},                                       \\
A_n & \defequal \sum_{1\le i\le n} \frac{1}{i},\; A_0\defequal 0, \quad
B_n \defequal \sum_{1\le i<j\le n}\frac{1}{ij},\; B_0=B_1\defequal 0,\quad 
C_n \defequal \sum_{1\le i\le n}\frac{1}{i^2}=A_n^2-2B_n. 
\end{align*}
We can write $S_{\Phi}^{-1}$ in terms of $I_0$, $I_1$, $I_2$, $I_3$. 
\[S_{\Phi}^{-1}=
  \begin{pmatrix}
  I_0 & \hbar\dot{I}_0 & \hbar I'_0 & \hbar^2\dot{I}'_0 \\
  I_1/\hbar &I_0+\dot{I}_1 & I'_1 & \hbar(I'_0+\dot{I}'_1) \\
  I_2/\hbar &\dot{I}_2 & I_0+I'_2 & \hbar(\dot{I}_0+\dot{I}'_2) \\
  I_3/\hbar^2 & (I_2+\dot{I}_3)\hbar & (I'_3+I_1+I_2)/\hbar &
  I_0+\dot{I}_1+\dot{I}_2+I'_2+\dot{I}'_3 
  \end{pmatrix},
\]
where $\dot{I}$ means $x\partial I/\partial x$ and 
$I'$ means $y\partial I/\partial y$.
We factorize $S_{\Phi}^{-1}$ as 
$S_{\Phi}^{-1}=S_{-}^{-1}S_{+}^{-1}$, where $S_{+}$
is a power series in $\hbar$ and $S_{-}$ is 
a power series in $\hbar^{-1}$ satisfying 
$S_{-}(\hbar=\infty)=\id$. 
We calculate $S_{\pm}^{-1}$ by Maple as
{\tiny 
\begin{align*}
 S_{+}^{-1}&=\begin{pmatrix} 
   1+2y+6y^2+20y^3+70y^4 & 2x(y+6y^2+30y^3+140y^4) & 
   2x(y+10y^2+70y^3+420y^4) & 
   8x^2(y^2+14y^3+126y^4)  \\
   0 & 1+2y+6y^2+20y^3+70y^4 & 2(y+7y^2+38y^3+187y^4) & 
   8x(y^2+11y^3+82y^4)     \\
   0 & 0 & 1+4y+16y^2+64y^3+256y^4  & 4x(y+8y^2+48y^3+256y^4) \\
   0 & 0 & 0 & 1+4y+16y^2+64y^3 
   \end{pmatrix}  \\
   &+\hbar 
   \begin{pmatrix} 
   0 & 0 & 2y+12y^2+60y^3+280y^4 & 2x(y+12y^2+90y^3+560y^4) \\
   0 & 0 & 0 & 2y+12y^2+60y^3+280y^4 \\
   0 & 0 & 0 & 0 \\
   0 & 0 & 0 & 0
   \end{pmatrix} 
   +O(y^5),  \\
S_{-}^{-1}&=
   \begin{pmatrix}
   1 & 0 & 0 & 0 \\
   0 & 1 & 0 & 0 \\
   0 & 0 & 1 & 0 \\
   0 & 0 & 0 & 1 
   \end{pmatrix} 
  +\frac{1}{\hbar}
   \begin{pmatrix}
   2x(y+4y^2+16y^3+64y^4) & 2x^2(y^2+8y^3+48y^4) & 2x^2(y^2+8y^3+48y^4)
   & 0 \\
   2y+5y^2+\frac{44}{3}y^3+\frac{93}{2}y^4 & x(-1+2y^2+12y^3+58y^4)
   & 2x(y^2+6y^3+29y^4) & 0 \\
   2y+3y^2+\frac{20}{3}y^3+\frac{35}{2}y^4 & x(1+2y+6y^2+20y^3+70y^4) 
   & 2x(y+3y^2+10y^3+35y^4) & 0 \\
   0 & 2y+3y^2+\frac{20}{3}y^3+\frac{35}{2}y^4 & 
   2y+3y^2+\frac{20}{3}y^3+\frac{35}{2}y^4 & 0
  \end{pmatrix} \\
  &+\frac{1}{\hbar^2} 
  \begin{pmatrix} 
  3x^2(y^2+8y^3+48y^4) & 4x^3(y^3+12y^4) & 4x^3(y^3+12y^4) & 0 \\
  x(-1-2y-2y^2+6y^3+\frac{190}{3}y^4) & x^2(\frac{1}{2}-2y^2-8y^3-15y^4) & 
  x^2(-2y^2-8y^3-15y^4) & 0 \\
  x(1+2y+10y^2+42y^3+\frac{514}{3}y^4) & x^2(-\frac{1}{2}+3y^2+24y^3+143y^4) & 
  x^2(3y^2+24y^3+143y^4) & 0 \\
  -2y+\frac{3}{2}y^2+\frac{88}{9}y^3+\frac{937}{24}y^4 & 
  x(-1-2y-2y^2+2y^3+\frac{94}{3}y^4) & x(-2y-2y^2+2y^3+\frac{94}{3}y^4) & 0 
  \end{pmatrix} \\
  &+\frac{1}{\hbar^3}
  \begin{pmatrix} 
  x^3(\frac{10}{3}y^3+40y^4) & 5x^4y^4 & 5x^4y^4 & 0 \\
  x^2(\frac{1}{4}-y-\frac{17}{2}y^2-42y^3-\frac{355}{2}y^4) &
  x^3(-\frac{1}{6}-y^2-14y^3-105y^4) & 
  x^3(-y^2-14y^3-105y^4) & 0 \\
  x^2(-\frac{1}{4}+y+\frac{15}{2}y^2+46y^3+\frac{499}{2}y^4) & 
  x^3(\frac{1}{6}+y^2+\frac{40}{3}y^3+113y^4) & 
  x^3(y^2+\frac{40}{3}y^3+113y^4) & 0 \\
  x(-2-4y-16y^2-53y^3-\frac{1552}{9}y^4) & 
  x^2(\frac{3}{4}-\frac{9}{2}y^2-32y^3-\frac{335}{2}y^4) & 
  x^2(-\frac{9}{2}y^2-32y^3-\frac{335}{2}y^4) & 0 
  \end{pmatrix}+O(\hbar^{-4},y^5).
\end{align*}
}
By the proof of Theorem \ref{thm:canonicalframe}, 
$S_{+}$ gives the gauge transformation $Q$ such that 
$\hPhi\defequal \Phi\circ Q$ is canonical. 
First, the connection matrices are transformed as follows:
\[\hOmega_1=S_{+}^{-1}\Omega_1S_{+}
       +\hbar x\frac{\partial S_{+}}{\partial x}, \quad 
  \hOmega_2=S_{+}^{-1}\Omega_1S_{+}
       +\hbar y\frac{\partial S_{+}}{\partial y}.
\]
We have 
{\tiny
\begin{eqnarray*} 
S_+&=&
\begin{pmatrix}
1-2y-2{y}^{2}-4{y}^{3}-10{y}^{4} &
-2xy-4x{y}^{2}-12x{y}^{3}-40x{y}^{4}& 
-2y\hbar-2xy-4x{y}^{2}-12x{y}^{3}-40x{y}^{4} &
-2xy\hbar \\
0 & 1-2y-2{y}^{2}-4{y}^{3}-10{y}^{4} & 
-2y-2{y}^{2}-4{y}^{3}-10{y}^{4} & 
-2y\hbar \\
0 & 0 & 1-4y & -4xy \\
0 & 0 & 0 & 1-4y 
\end{pmatrix}+O(y^5), 
\\
\hOmega_1&=&
\begin{pmatrix} 
2xy+8x{y}^{2}+32x{y}^{3}+128x{y}^{4} &
4{x}^{2}{y}^{2}+32{x}^{2}{y}^{3}+192{x}^{2}{y}^{4}&
4{x}^{2}{y}^{2}+32{x}^{2}{y}^{3}+192{x}^{2}{y}^{4}&
0  &  \\
1  & -x+2x{y}^{2}+12x{y}^{3}+58x{y}^{4}& 
2x{y}^{2}+12x{y}^{3}+58x{y}^{4} & 0  \\ 
0 & x+2xy+6x{y}^{2}+20x{y}^{3}+70x{y}^{4}& 
2xy+6x{y}^{2}+20x{y}^{3}+70x{y}^{4} & 0 \\ 
0 & 0 & 1 & 0 
\end{pmatrix}+O(y^5),
\\
\hOmega_2&=&
\begin{pmatrix} 
2xy+16x{y}^{2}+96x{y}^{3}+512x{y}^{4} & 
4{x}^{2}{y}^{2}+48{x}^{2}{y}^{3}+384{x}^{2}{y}^{4} & 
4{x}^{2}{y}^{2}+48{x}^{2}{y}^{3}+384{x}^{2}{y}^{4} & 
0 \\
2y+10{y}^{2}+44{y}^{3}+186{y}^{4} & 
4x{y}^{2}+36x{y}^{3}+232x{y}^{4} & 
4x{y}^{2}+36x{y}^{3}+232x{y}^{4} & 0  \\
1+2y+6{y}^{2}+20{y}^{3}+70{y}^{4} & 
2xy+12x{y}^{2}+60x{y}^{3}+280x{y}^{4} & 
2xy+12x{y}^{2}+60x{y}^{3}+280x{y}^{4} & 0 \\
0 & 1+2y+6{y}^{2}+20{y}^{3}+70{y}^{4} & 
1+2y+6{y}^{2}+20{y}^{3}+70{y}^{4} & 0 
\end{pmatrix}+O(y^5). 
\end{eqnarray*}
}
Second, we must transform coordinates.  
We can read the coordinate transformation (mirror transformation) 
from $S_{+}$ and $\hOmega_a$. 
The mirror transformation is determined by 
\[ \hOmega_a(1)=-\lparfrac{F}{q^a}+
   \sum_{b=1}^r\lparfrac{\log \hq^b}{q^a}p_b.
\]
(See equations (\ref{eq:mirrortransmethod1}),
 (\ref{eq:mirrortransmethod2}).)
From this, we calculate the transformation as  
\begin{gather} 
\label{eq:mirrortrans}
   x=\hx(1-2\hy+\hy^2), \quad  
   y=\hy-2\hy^2+3\hy^3-4\hy^4+5\hy^5+\cdots, 
\end{gather}
where $\hx=\hq^1$ and $\hy=\hq^2$. 
The function $F$ can be calculated as 
\begin{gather}
\label{eq:Fexample}
   F=-2x(y+4y^2+16y^3+64y^4+\cdots)=-2\hx\hy.  
\end{gather} 
The connection matrices are further transformed as 
\[ \hOmega_{\hat{1}}=\lparfrac{\log x}{\hx}\hOmega_1+
           \lparfrac{\log y}{\hx}\hOmega_2+\lparfrac{F}{\hx}, 
   \quad 
   \hOmega_{\hat{2}}=\lparfrac{\log x}{\hy}\hOmega_1+
           \lparfrac{\log y}{\hy}\hOmega_2+\lparfrac{F}{\hy}.
\]
We have 
\[
\hOmega_{\hat{1}}=
\begin{pmatrix}
0 & 4\hx^2\hy^2   & 4\hx^2\hy^2 & 0 \\
1 & -\hx(1-\hy^2) & 2\hx\hy^2   & 0 \\
0 & \hx(1-\hy^2)  & -2\hx\hy^2  & 0 \\
0 & 0             & 1           & -2\hx\hy
\end{pmatrix}, \;\,
\hOmega_{\hat{2}}=
\begin{pmatrix}
0 & 4\hx^2\hy^2 & 4\hx^2\hy^2 & 0 \\
0 & 2\hx\hy^2   & 4\hx\hy^2   & 0 \\
1 & -2\hx\hy^2  & -4\hx\hy^2  & 0 \\
0 &  1          & 1-2(\hy+\hy^2+\hy^3+\cdots) & -2\hx\hy  
\end{pmatrix}.  
\]
They are $\hbar$-independent, therefore 
the frame $\hPhi=\Phi\circ S_{+}$ is canonical. 
The function $f$ in the equation (\ref{eq:transformationofJ})
is given by $f=S_{+}(1)=1-2y-2y^2-4y^3-10y^4-\cdots$. 

In the above calculation, we solve for each 
coefficient of $x^iy^j$ recursively. 
Thus we do not reach the exact solution by this method.  
However, in this case, we can obtain the exact solution 
by an analytical method. 
The $(1/\hbar)$-expansion of $J_{\Phi}$ is 
\begin{equation*}
J_{\Phi}=e^{(p_1\log x+p_2\log y)/\hbar}
       \left(\frac{1}{f}+\frac{1}{\hbar}(H_1p_1+H_2p_2-\frac{F}{f})
       +o(\hbar^{-1})\right), 
\end{equation*} 
where 
\begin{align*}
  f(x,y)\defequal& \left(\sum_{n=0}^\infty 
                        \frac{(2n)!}{(n!)^2}y^n\right)^{-1}
  =\sqrt{1-4y}, 
  \quad 
  F(x,y)\defequal xy\frac{\partial f/\partial y}{f}
  =-\frac{2xy}{1-4y}, \\
  H_1(x,y)\defequal& \sum_{n=0}^\infty \frac{(2n)!}{(n!)^2}A_ny^n
  =\frac{2}{\sqrt{1-4y}}\log (\frac{1}{2}+\frac{1}{2\sqrt{1-4y}}), \\
  H_2(x,y)\defequal& 2\sum_{n=0}^\infty \frac{(2n)!}{(n!)^2}
                    (A_{2n}-A_n)y^n
          =-\frac{2}{\sqrt{1-4y}}\log(\frac{1+\sqrt{1-4y}}{2}). 
\end{align*}
This asymptotics and Proposition 
\ref{prop:giventalmirrortrans} show that the mirror transformation
is given by 
\[J_{\Phi}\longmapsto fe^{F/\hbar}J_{\Phi},\quad 
\log \hx=\log x+fH_1,\quad  \log\hy=\log y+fH_2.
\] 
The inverse change of variables is given by the 
simple formulas. 
\[x=\hx(1-\hy)^2,\quad y=\frac{\hy}{(1+\hy)^2}.
\] 
From this, it follows that $f=(1-\hy)/(1+\hy)$ and $F=-2\hx\hy$. 
They coincide with the equations (\ref{eq:mirrortrans}), 
(\ref{eq:Fexample}).  
We transform $\Delta^{\mV}$ and $P_a$ as 
\begin{align*}
 \hat{\Delta}^{\mV}=&f\Delta^{\mV}, \\
 \hat{P}_1=&\lparfrac{\log x}{\hx}P_1+
           \lparfrac{\log y}{\hx}P_2+\lparfrac{F}{\hx}
          =P_1-2\hx\hy,  \\
 \hat{P}_2=&\lparfrac{\log x}{\hy}P_1+
           \lparfrac{\log y}{\hy}P_2+\lparfrac{F}{\hy}
          =-\frac{2\hy}{1-\hy}P_1+\frac{1-\hy}{1+\hy}P_2-2\hx\hy.
\end{align*} 
Then the Picard-Fuchs equations (\ref{eq:Picard-Fuchsofsuperspace}) 
is transformed as 
\begin{gather*} 
{\hat{P}_1}^2\hat{\Delta}^{\mV}=4\hx^2\hy^2\hat{\Delta}^{\mV}
  -\hx(1-\hy^2)(\hat{P}_1\hat{\Delta}^{\mV}
                 -\hat{P}_2\hat{\Delta}^{\mV}), \\
2\hx\hy(1+\hy)(\hat{P}_1\hat{\Delta}^{\mV}-\hat{P}_2\hat{\Delta}^{\mV})
  +(1-3\hy)\hat{P}_1\hat{P}_2\hat{\Delta}^{\mV} 
  +2\hy\hat{P}_1^2\hat{\Delta}^{\mV}
  -(1-\hy)\hat{P}_2^2\hat{\Delta}^{\mV}=0, 
\end{gather*}  
where, by the abuse of notation, 
$\hx$, $\hy$ mean the pull-backs by
the corresponding covering transformations. 
We define a new frame $\hPhi$ as 
\begin{gather*}
  \hPhi(1)=\hat{\Delta}^{\mV}, \;\; 
  \hPhi(p_1)=\hat{P}_1\hat{\Delta}^{\mV}, \;\; 
  \hPhi(p_2)=\hat{P}_2\hat{\Delta}^{\mV}, \\ 
  \hPhi(p_1p_2)=\hat{P}_1\hat{P}_2\hat{\Delta}^{\mV}
                +2\hx\hy^2(\hat{P}_2-\hat{P}_1)\hat{\Delta}^{\mV}
                -4\hx^2\hy^2\hat{\Delta}^{\mV}. 
\end{gather*}
By using the Picard-Fuchs equations, 
we can obtain the same connection matrices $\hOmega_{\hat{a}}$ 
as above, but the $(4,3)$-entry $1-2(\hy+\hy^2+\hy^3+\cdots)$ of 
$\hOmega_{\hat{2}}$ 
is replaced by the exact formula $(1-3\hy)/(1-\hy)$.

\end{document}